\documentclass{amsart}
\usepackage{graphicx} 
\usepackage[english]{babel}
\usepackage[margin=2.5cm]{geometry} 
\usepackage{amssymb, amsfonts, amsthm, amsmath}
\usepackage[hidelinks]{hyperref}
\usepackage{caption} 
\usepackage[dvipsnames]{xcolor} 
\usepackage{thmtools,tikz-cd}

\usepackage{subcaption}

\usepackage{bm}
\usepackage{pgf,tikz} 
\usetikzlibrary{matrix,calc,positioning,intersections,through,angles,patterns}
\usetikzlibrary{arrows.meta}

\usetikzlibrary{math,arrows,through, quotes}
\usepackage[mathscr]{euscript}
\usepackage{enumitem} 

\usepackage{caption} 

\usepackage{parskip}
\usepackage{enumitem}

\usepackage{changepage} 
\newenvironment{statement}
{	\itshape
	\vspace{1mm}
	\begin{adjustwidth}{20mm}{20mm}
		\begin{center}
			\large
}
{ 
		\end{center}
	\end{adjustwidth}
	\vspace{1mm}
}

\usepackage{calc}

\usepackage{ifthen}

\usepackage{bm} 

\usepackage{cancel} 

\usepackage{thmtools}
\declaretheorem[style=plain,name=Theorem,qed={\tiny$\blacksquare$},numberwithin=section]{theorem}
\declaretheorem[style=plain,name=Corollary,sibling=theorem,qed={\tiny$\blacksquare$}]{corollary}
\declaretheorem[style=definition,name=Definition,sibling=theorem,qed={\tiny$\blacksquare$}]{definition}
\declaretheorem[style=plain,name=Lemma,sibling=theorem,qed={\tiny$\blacksquare$}]{lemma}
\declaretheorem[style=definition,name=Example,sibling=theorem,qed={\tiny$\blacksquare$}]{example}
\declaretheorem[style=definition,name=Remark,sibling=theorem,qed={\tiny$\blacksquare$}]{remark}

\declaretheorem[style=plain,name=Proposition,sibling=theorem,qed={\tiny$\blacksquare$}]{proposition}

\numberwithin{equation}{section}

\usepackage[inline]{trackchanges} 
\addeditor{TG}
\addeditor{SR}
\addeditor{NA}
\addeditor{MG}

\newcommand{\Span}{\operatorname{span}}

\DeclareMathSymbol{\shortminus}{\mathbin}{AMSa}{"39}

\DeclareMathOperator{\cro}{cr} 
\DeclareMathOperator{\sr}{sr} 
\DeclareMathOperator{\mr}{mr} 

\newcommand{\tcd}{\mathcal T} 
\newcommand{\lrp}{{\mathcal{L}\kern-0.15em\mathcal{R}}} 
\newcommand{\pb}{{G}} 
\newcommand{\edgeweight}{\mu} 
\newcommand{\enm}[2]{{\mathcal S_{#2}^{#1}}} 

\newcommand{\tcdm}{{\mathcal T_{\bm{\circlearrowright}}}} 
\newcommand{\tcdp}{{\mathcal T_{\bm{\circlearrowleft}}}} 
\newcommand{\vrc}{{\mathsf{v}}} 
\newcommand{\relation}{{\mathsf{R}}} 
\newcommand{\point}{{\mathsf{P}}} 
\newcommand{\dimension}{d} 
\newcommand{\hyperplane}{{\mathsf{H}}} 
\newcommand{\projectiveline}{{\mathsf{L}}}
\newcommand{\centerofprojection}{{\mathsf{Z}}}
\newcommand{\affinequiver}{\quiver^{\mathrm{aff}}}
\newcommand{\projectivequiver}{\quiver^{\mathrm{pro}}}
\newcommand{\pbint}{{\pb^{\mathrm{int}}}}

\newcommand{\disk}{\mathbb D}
\newcommand{\complex}{\mathfrak{F}}
\newcommand{\ncp}{\pi}

\DeclareMathOperator{\pro }{\mathbf{Pro}} 
\DeclareMathOperator{\aff }{\mathbf{Aff}} 

\newcommand{\an}{{\bm A}} 
\newcommand{\ban}{{\bm B}} 

\newcommand{\rank}{{rank}}
\newcommand{\mrank}{\operatorname{rk}}
\newcommand{\maxrank}{\mrank}

\newcommand{\plabic}{{BTB}}
\newcommand{\Plabic}{{BTB}}
\newcommand{\stargraph}{{\star(\pb)}} 
\newcommand{\starb}{b^\star} 
\newcommand{\starw}{w^\star} 
\newcommand{\starf}{f^\star} 
\newcommand{\starvrc}{\vrc^\star} 
\newcommand{\starrelation}{\relation^\star} 

\newcommand{\tcdmap}{\mathsf{T}}
\newcommand{\quiver}{Q}
\newcommand{\maps}{\mathrm{Moduli}}
\newcommand{\subspacea}{\mathsf{U}}
\newcommand{\subspaceb}{\mathsf{V}}

\newcommand{\Z}{\mathbb{Z}} 
\newcommand{\R}{\mathbb{R}} 
\newcommand{\C}{\mathbb{C}} 

\newcommand{\p}[1]{\bm{#1}}

\newcommand{\CP}{\mathbb{C}\mathsf{P}} 
\newcommand{\resplitmove}{{resplit}}
\newcommand{\resplitmovecap}{{Resplit}}

\newcommand{\menge}[2][]{ 
	  \if\relax\detokenize{#1}\relax
		  \left\{#2\right\}
	  \else
		  \left\{#1\ | \ #2\right\}
	  \fi}

\def\centerarc[#1](#2)(#3:#4:#5)
{ \draw[#1] ($(#2)+({#5*cos(#3)},{#5*sin(#3)})$) arc (#3:#4:#5); }

\tikzset{rvert/.style={draw,circle,fill=red,minimum size=5pt,inner sep=0pt} }
\tikzset{redge/.style={draw,densely dashed, red} }
\tikzset{wqgvert/.style={draw,black,rectangle,fill=white,minimum size=4.5pt,inner sep=0pt}  } 
\tikzset{bqgvert/.style={draw,black,rectangle,fill=black,minimum size=4.5pt,inner sep=0pt}  } 
\tikzset{qvert/.style={draw,black,circle,fill=gray,minimum size=5pt,inner sep=0pt}  } 
\tikzset{bvert/.style={draw,circle,fill=black,minimum size=5pt,inner sep=0pt}  } 
\tikzset{wvert/.style={draw,circle,fill=white,minimum size=5pt,inner sep=0pt}  } 
\tikzset{nvert/.style={}  } 
\tikzset{hvert/.style={draw,circle,minimum size=10pt,inner sep=0pt,line width=0.8pt,gray!75}  } 
\tikzset{hqvert/.style={draw,black,circle,pattern color=gray,pattern=north west lines,minimum size=7pt,inner sep=0pt}  } 
\tikzset{wfoc/.style={draw,rectangle,fill=white,minimum size=4pt,inner sep=0pt}  } 
\tikzset{bfoc/.style={draw,rectangle,fill=black,minimum size=4pt,inner sep=0pt}  } 
\tikzset{gfoc/.style={draw,rectangle,fill=gray,minimum size=4pt,inner sep=0pt}  } 

\newcommand{\tztriangle}[3]{
	\begin{tikzpicture}[baseline={([yshift=-.7ex]current bounding box.center)},scale=.8] 
		\tikzstyle{bvert}=[draw,circle,fill=black,minimum size=2.5pt,inner sep=0pt]
		
		\node[bvert] (v1)  at (0,0) {};
		\node[bvert] (v2) at (1,0) {};
		\node[bvert] (v3) at (.5,0.875) {};
		
		\ifnum#1=1 \draw[-] (v2) -- (v3);\else \draw[-,gray!50] (v2) -- (v3); \fi
		\ifnum#2=1 \draw[-] (v1) -- (v3);\else \draw[-,gray!50] (v1) -- (v3); \fi
		\ifnum#3=1 \draw[-] (v1) -- (v2);\else \draw[-,gray!50] (v1) -- (v2); \fi
\end{tikzpicture}}

\newcommand{\tzstar}[3]{
	\begin{tikzpicture}[baseline={([yshift=-.7ex]current bounding box.center)},scale=.8] 
		\tikzstyle{bvert}=[draw,circle,fill=black,minimum size=2.5pt,inner sep=0pt]
		
		\node[bvert] (v1)  at (0,0) {};
		\node[bvert] (v2) at (1,0) {};
		\node[bvert] (v3) at (.5,0.875) {};
		\node[bvert] (v4) at (.5,0.35) {};
		
		\ifnum#1=1 \draw[-] (v4) -- (v1); \else \draw[-,gray!50] (v4) -- (v1); \fi
		\ifnum#2=1 \draw[-] (v4) -- (v2); \else \draw[-,gray!50] (v4) -- (v2); \fi
		\ifnum#3=1 \draw[-] (v4) -- (v3); \else \draw[-,gray!50] (v4) -- (v3); \fi
		
\end{tikzpicture}}

\newcommand{\tzspider}[3][0]{
	\begin{tikzpicture}[baseline={([yshift=-.7ex]current bounding box.center)},scale=1,rotate=#1] 
		\tikzstyle{bvert}=[draw,circle,fill=black,minimum size=2.5pt,inner sep=0pt]
		
		\node[bvert] (v1)  at (0,0) {};
		\node[bvert] (v2) at (1,0) {};
		\node[bvert] (v3) at (1,1) {};
		\node[bvert] (v4) at (0,1) {};
		\node[bvert] (w1) at (.25,.25) {};
		\node[bvert] (w3) at (.75,.75) {};
		
		\ifnum#2=1 \draw[-] (w1) -- (v1);\else \draw[-,gray!50] (w1) -- (v1); \fi
		\ifnum#2=2 \draw[-] (w1) -- (v2);\else \draw[-,gray!50] (w1) -- (v2);  \fi
		\ifnum#2=3 \draw[-] (w1) -- (v4);\else \draw[-,gray!50] (w1) -- (v4);  \fi
		\ifnum#3=1 \draw[-] (w3) -- (v3);\else \draw[-,gray!50] (w3) -- (v3);  \fi
		\ifnum#3=2 \draw[-] (w3) -- (v4);\else \draw[-,gray!50] (w3) -- (v4);  \fi
		\ifnum#3=3 \draw[-] (w3) -- (v2);\else \draw[-,gray!50] (w3) -- (v2);  \fi
\end{tikzpicture}}

\def\centerarc[#1](#2)(#3:#4:#5)
{ \draw[#1] ($(#2)+({#5*cos(#3)},{#5*sin(#3)})$) arc (#3:#4:#5); }

\usetikzlibrary{decorations.pathreplacing,decorations.markings}
\tikzset{ 
	mid arrow/.style={postaction={decorate,decoration={
				markings,
				mark=at position  0.5*\pgfdecoratedpathlength+2.5pt with {\arrow[line width=0.4 mm]{latex}}
	}}},
	mid rarrow/.style={postaction={decorate,decoration={
				markings,
				mark=at position 0.5*\pgfdecoratedpathlength+2.5pt with {\arrow[line width=0.4 mm]{latex reversed}}
	}}},
}
\tikzset{ 
	orient/.style={thick,postaction={decorate,decoration={
				markings,
				mark=at position 0.55*\pgfdecoratedpathlength+2pt with {\arrow{Triangle}}
	}}},
	rorient/.style={thick,postaction={decorate,decoration={
				markings,
				mark=at position 0.45*\pgfdecoratedpathlength+1pt with {\arrow{Triangle[reversed]}}
	}}},
}
\tikzset{ 
	orientb/.style={thick,postaction={decorate,decoration={
				markings,
				mark=at position 0.5*\pgfdecoratedpathlength+1pt with {\arrow{Classical TikZ Rightarrow}},
				mark=at position 0.5*\pgfdecoratedpathlength-1pt with {\arrow{Classical TikZ Rightarrow}},				
	}}},
	rorientb/.style={thick,postaction={decorate,decoration={
				markings,
				mark=at position 0.5*\pgfdecoratedpathlength+1pt with {\arrow{Classical TikZ Rightarrow[reversed]}},
				mark=at position 0.5*\pgfdecoratedpathlength-1pt with {\arrow{Classical TikZ Rightarrow[reversed]}}				
	}}},
}

\tikzset{circle through 3 points/.style n args={3}{%
		insert path={let    \p1=($(#1)!0.5!(#2)$),
			\p2=($(#1)!0.5!(#3)$),
			\p3=($(#1)!0.5!(#2)!1!-90:(#2)$),
			\p4=($(#1)!0.5!(#3)!1!90:(#3)$),
			\p5=(intersection of \p1--\p3 and \p2--\p4)
			in },
		at={(\p5)},
		circle through= {(#1)}
}}

\title[Multiple cluster algebra structures for TCD maps I: theoretical framework]{Multiple cluster algebra structures for triple crossing diagram maps I: theoretical framework}

\date{\today}

\author{Niklas Affolter}
\address{Technische Universit\"at Wien, Institute of Discrete Mathematics and Geometry, Wiedener Hauptstr. 8-10/104, A-1040 Vienna,
Austria}
\email{affolter@posteo.net}
\author{Terrence George}
\address{Department of Mathematics, Massachusetts Institute of Technology, 77 Massachusetts Ave., Cambridge, MA
02139, USA}
\email{tegeorge@mit.edu}
\author{Max Glick}
\address{Google Inc, Pittsburgh, PA 15206, USA}
\email{max.i.glick@gmail.com}
\author{Sanjay Ramassamy}
\address{Universit\'e Paris-Saclay, CNRS, CEA, Institut de Physique Th\'eorique, 91191 Gif-sur-Yvette, France}
\email{sanjay.ramassamy@ipht.fr}

\begin{document}

\begin{abstract}
We introduce triple crossing diagram (TCD) maps, which encode projective configurations of points and lines, as a unified framework for constructions arising in various areas of geometry, such as discrete differential geometry, discrete geometric dynamics and hyperbolic geometry. We define two types of local moves for TCD maps, one of which is governed by the discrete Schwarzian KP (dSKP) equation, and establish their multi-dimensional consistency. We construct two distinct cluster structures on the space of TCD maps, called projective and affine cluster structures, and show that they are related via an operation called section. This framework organizes and unifies a wide range of examples, including Q-nets, Darboux maps, line complexes, T-graphs, t-embeddings, triangulations and geometric discrete integrable systems such as the pentagram map and cross-ratio dynamics, which are further developed in the companion paper~\cite{paper2} and in~\cite{agrcrdyn}.
\end{abstract}

\maketitle

\setcounter{tocdepth}{1}
\tableofcontents

\newpage

\section{Introduction}
\label{sec:intro}

\emph{Cluster algebras}, introduced by Fomin and Zelevinsky~\cite{MR1887642}, provide a unifying framework for various algebraic and geometric objects in terms of \emph{clusters} of variables related by explicit local transformations called \emph{mutations}. The underlying combinatorial data of a cluster algebra is sometimes called a \emph{cluster structure}. Concretely, a cluster structure consists of an oriented graph, called a \emph{quiver}, whose vertices carry \emph{cluster variables}. In this paper, we focus on quivers arising as oriented duals of planar bipartite graphs and on cluster variables of coefficient type, called $y$-variables by~\cite{fzcoefficients} and $X$-variables by \cite{fghighertm}. Mutation acts locally on the quiver and its variables, producing a new quiver and cluster variables according to certain fixed combinatorial rules.

Let us mention three major directions of research that have involved such cluster structures:
\begin{itemize}
\item \textbf{Grassmannians:} the parametrization of the totally positive Grassmannian and its positroid cells via plabic graphs~\cite{postgrass,scottgrass, GalashinLamcluster}. 
\item \textbf{Discrete dynamics:} the integrability of certain discrete-time dynamical systems, such as the {pentagram map} \cite{schwartz,ostpentagram,glickpentagram} and its generalizations, {cross-ratio dynamics} \cite{afitcrossratio, agrcrdyn}, polygon recutting \cite{adler_recutting,recutting}, etc.
\item \textbf{Dimer embeddings:} planar realizations of the dimer model and related embedding theories, including {T-graphs} \cite{kenyonsheffield}, {t-embeddings} \cite{amiquel,kenyonlam,clrtembeddings}, {s-embeddings} \cite{chelkaksembeddings}, {S-graphs} \cite{chelkaksgraphs} and {Tutte embeddings} \cite{tutteembedding}.
\end{itemize}

Recently, the notion of a \emph{vector-relation configuration (VRC)} was introduced in~\cite{vrc}. VRCs possess a cluster strcuture and provide a unified framework encompassing both the Grassmannian and discrete dynamics cases. However, the cluster structures arising from dimer embeddings exhibit fundamentally different behavior and do not fit in the VRC framework. The purpose of this paper is to introduce a modification of VRCs, which we call \emph{triple crossing diagram maps (TCD maps)}. We show that TCD maps naturally possess a second cluster structure that captures the dimer embeddings, thus providing a framework with wider scope. Moreover, we show that the second cluster structure is related to the first cluster structure via a geometric operation called \emph{section}. Finally, we prove the \emph{multi-dimensional consistency} of TCD maps, a notion of integrability from discrete differential geometry~\cite{absintquads,absquads, ddgbook}, which reveals their close connection to \textbf{discrete integrable systems}, in particular to Desargues maps~\cite{doliwadesargues,doliwadesarguesweyl} and dSKP lattices~\cite{dndskp,ncwqdskp,ksclifford}. Note that here ``discrete integrability'' is a different notion of integrability from the discrete dynamics mentioned above, which are integrable in the sense of Liouville-Arnold.
In this way, TCD maps also encompass a fourth direction of research, thereby also establishing the existence of cluster structures in these systems.

\subsection{TCD maps and their projective cluster structures}

\begin{figure}
    \centering
    \includegraphics[angle=90,scale=0.95]{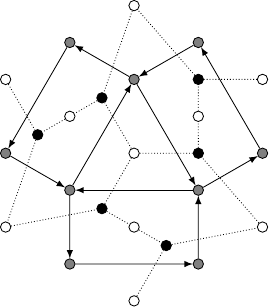}
    \hspace{5mm}
    \includegraphics[angle=90,scale=0.95]{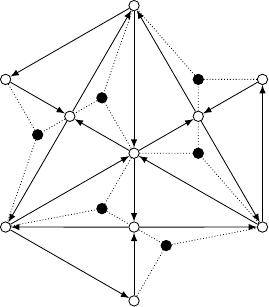}
    \hspace{5mm}
    \includegraphics[angle=90,scale=0.95]{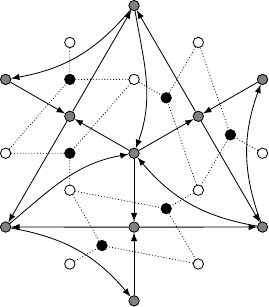}
    \caption{A BTB graph $\pb$ (dotted) with its projective quiver (left), its affine quiver (center), and the projective quiver of a section $\sigma(\pb)$ (right). }
    \label{fig:quivers}
\end{figure}

A \emph{BTB graph} is a black-trivalent bipartite graph, namely a planar bipartite graph $\pb$ such that every black vertex has degree 3. For the sake of exposition and since most examples of interest fall into this class, we assume unless otherwise stated that all graphs considered are \emph{minimal} in the sense of \cite{thurstontriple, postgrass}, although many of our definitions below remain valid without this assumption. A \emph{TCD map} is a map from the white vertices of $\pb$ to a fixed projective space $\CP^\dimension$ such that, at every black vertex, the images of the three adjacent white vertices lie on a common projective line. The name reflects that BTB graphs are in bijection with Thurston's triple crossing diagrams (TCDs)~\cite{thurstontriple}.

\begin{figure}
    \centering

    \raisebox{-0.45\height}{\includegraphics[scale=1]{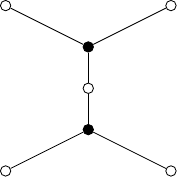}}
    \hspace{1mm}$\leftrightarrow$\hspace{1mm}
    \raisebox{-0.45\height}{\includegraphics[scale=1,angle=90]{tikz/introresplit.pdf}}
    \hspace{15mm}
    \raisebox{-0.45\height}{\includegraphics[scale=1]{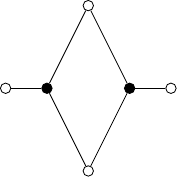}}
    \hspace{1mm}$\leftrightarrow$\hspace{1mm}
    \raisebox{-0.45\height}{\includegraphics[scale=1,angle=90]{tikz/spider1.pdf}}

    \caption{Local moves for BTB graphs: resplit (left) and spider move (right). 
    }
    \label{fig:introlocalmoves}
\end{figure}

TCD maps are the projectivized form of VRCs together with the additional trivalence condition at black vertices (whereas VRCs are defined on general planar bipartite graphs called \emph{plabic graphs}). VRCs have several local moves: the \emph{spider move} and the mutually inverse operations of \emph{splitting} and \emph{contracting} vertices. When restricted to BTB graphs, only two moves remain: the \emph{spider move} and the \emph{resplit}, which is the composition of a contraction followed by a split (see~Figure~\ref{fig:introlocalmoves}). Moreover, any plabic graph can be transformed into a BTB graph through a sequence of splits and contractions; thus, BTB graphs, although being a subclass, possess the same scope as plabic graphs. On the level of TCD maps, the spider move leaves all white vertices fixed, acting merely as a reparametrization of the TCD map, while the resplit replaces a white vertex with a new one which is mapped to the intersection of the two lines associated with the adjacent black vertices.

\begin{figure}
    \centering
    \includegraphics[scale=0.95]{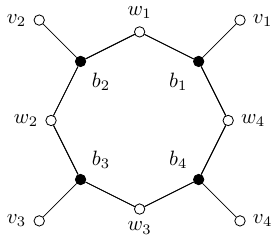}
    \hspace{20mm}
    \includegraphics[scale=0.95]{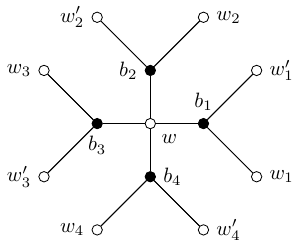}
    \caption{Local configurations in the definition of the projective (left) and affine (right) cluster variables.}
    \label{fig:localconfig}
\end{figure}

As mentioned earlier, TCD maps possess two cluster structures. The first cluster structure, which we call the \emph{projective cluster structure}, is the cluster structure for VRCs was introduced in~\cite{vrc} reinterpreted in the framework of TCD maps. The \emph{projective quiver} of a BTB graph $\pb$ is the oriented dual graph of $\pb$, so the cluster variables are naturally assigned to the faces of $\pb$ (see~Figure~\ref{fig:quivers}~(left)). For a face~$f$, the corresponding \emph{projective cluster variable}~$X_f$ is defined in terms of a {multi-ratio}:
\[
    X_f   := -\prod_{i=1}^m \frac{\tcdmap(w_{i-1}) - \tcdmap(v_{i})}{\tcdmap(w_{i}) - \tcdmap(v_{i})},
\]
with the labeling of vertices as in Figure~\ref{fig:localconfig}~(left) and the expression on the right is computed in any affine chart of~$\CP^\dimension$. The terminology “projective’’ reflects that the cluster variables, being multi-ratios, are invariant under projective transformations. It was shown in~\cite{vrc} that the projective cluster variables transform under the spider move according to the mutation formula for $X$ cluster variables and that they are invariant under resplits.

\subsection{Affine cluster structures}

Our first contribution is the introduction of a second cluster structure for TCD maps, which we call the \emph{affine cluster structure}. Fix a hyperplane $\hyperplane \subset \CP^\dimension$ and consider a TCD map $\tcdmap$ in the affine chart of $\CP^\dimension$ for which $\hyperplane$ is the hyperplane at infinity. The \emph{affine quiver} has a vertex for each white vertex of $\pb$ (see~Figure~\ref{fig:quivers}~(center)) and the associated \emph{affine cluster variable} is defined by the following \emph{star-ratio}:
\[
    Y_w := - \prod_{i=1}^m \frac{\tcdmap(w_i)-\tcdmap(w)}{\tcdmap(w'_i)-\tcdmap(w)}, 
\]
with the labeling of vertices as in Figure~\ref{fig:localconfig}~(right). The terminology “affine’’ reflects that these cluster variables are invariant under affine transformations (i.e. projective transformations that fix $\hyperplane$). We show that affine cluster variables transform under the resplit according to the mutation formula for $X$ cluster variables but are invariant under {spider moves}. In this sense, affine and projective cluster structures exhibit opposite behavior: the affine structure mutates under resplits but not spider moves, whereas the projective structure mutates under spider moves but not resplits. Since the definition of affine cluster variables requires black vertices to be trivalent, this provides an additional reason for working with BTB graphs rather than general plabic graphs. 

The affine cluster structure is what brings embeddings, Grassmannians and discrete dynamics under a unified framework: in particular, we show that any t-embedding (a special case of a dimer embedding) can be realized as a TCD map $\tcdmap$ whose known cluster structure coincides with the affine cluster structure of $\tcdmap$.

\subsection{Sections and the correspondence between cluster structures}

Our second contribution is the introduction of a geometric operation on TCD maps called a \emph{section}. As before, let $\hyperplane \subset \CP^\dimension$ be a hyperplane. Informally, a section is a map $\sigma_\hyperplane(\tcdmap)$ that assigns to each black vertex~$b$ of~$\pb$ the point where the line associated with~$b$ meets~$\hyperplane$. We show that a section is again a TCD map on one of several possible related graphs $\sigma(\pb)$, where the different choices of $\sigma(\pb)$ are related by resplits. Closely related constructions also appear in~\cite{galashincritical, Galashin2024origami}. There, it is shown that two iterations of the section preserves positivity of the cluster variables. The combinatorial operation sending $\pb$ to $\sigma(\pb)$ is also called the \emph{shift-by-1 map} \cite{Galplabictilings} and \emph{$T$-duality} \cite{MR4651900, CasalsLeShermanBennettWeng2023}.

If $\tcdmap$ is a map to $\CP^\dimension$, we refer to $\dimension$ as its \emph{rank}. Since a section of $\tcdmap$ has rank $\dimension-1$, we can iteratively take $\dimension-2$ sections, each of which possesses a projective and an affine cluster structure. Our third contribution is the following correspondence between projective and affine cluster structures:

\begin{statement}
    The affine cluster structure of a TCD map agrees with the projective cluster structure of a section. 
\end{statement}

See also Figure~\ref{fig:quivers}~(center and right) and Theorem~\ref{th:affprojcluster}. Therefore, TCD maps have $\dimension$ cluster structures and every affine cluster structure is also a projective cluster structure though of a different TCD map, except in the case of TCD maps of rank 1, where no nontrivial sections exist; notably, this case includes the important example of t-embeddings.

\subsection{Multi-dimensional consistency}

Given a BTB graph~$\pb$, its \emph{move graph} is the graph whose vertices represent all BTB graphs obtained from~$\pb$ by sequences of local moves, with edges corresponding to individual moves. Cycles in the move graph based at~$\pb$ thus correspond to sequences of moves that return~$\pb$ to itself. Thurston~\cite{thurstontriple} conjectured that all such cycles are generated by three elementary ones, a result recently proved by Balitsky and Wellmann~\cite{bwtriple}. Our fourth main contribution is that TCD maps exhibit the following property:

\begin{statement}\label{thm:intro_consistency}
	TCD maps are invariant under any sequence of moves forming a cycle.
\end{statement}

The proof is a consequence of the classical {Desargues' theorem}. This property is known as \emph{multi-dimensional consistency}; see Theorem~\ref{th:tcd_consistency_along_fund_cycles}. 

Typically, in discrete integrable systems multi-dimensional consistency is studied on maps defined on lattices. There is a well-known way to label vertices of $\pb$, introduced in the context of the cluster structure on the Grassmannian~\cite{scottgrass} and studied in relation to weak separation~\cite{OPS}. In our context, this labeling realizes TCD maps as maps on the root lattice of type $A_n$. This lattice perspective allows us to compare TCD maps with other geometric maps defined on $A_n$ lattices studied in discrete integrable systems, which results in our fifth contribution:

\begin{statement}
	TCD maps are the restrictions of \emph{Desargues maps}. Moreover, coordinate projections of TCD maps form \emph{dSKP lattices}.
\end{statement}

See also Theorem~\ref{th:tcddesargues} and Theorem~\ref{th:tcddskp}, where we show that this statement is also compatible with local moves in an appropriate sense.
Desargues maps were introduced by Doliwa~\cite{doliwadesargues, doliwadesarguesweyl} as maps to~$\CP^\dimension$ defined by incidence relations, while dSKP lattices arise as a discretization of the Schwarzian KP equation~\cite{dndskp,ncwqdskp,ksclifford}.

\subsection{Applications and related work}

In the companion paper~\cite{paper2}, we present many geometric systems that arise as special cases of TCD maps. This list contains many classical integrable examples from discrete differential geometry, including Q-nets, Darboux maps, line complexes, line compounds, Laplace transformations, the pentagram map and its generalizations, ideal hyperbolic triangulations, projective flag configurations, t-embeddings (see also Section~\ref{sec:tembeddings}) and T-graphs. Some of these examples have already appeared in~\cite{vrc}. 

Much of the material in the present paper and its companion paper~\cite{paper2} appeared in the thesis of the first author \cite{affolterthesis}. The article~\cite{agrcrdyn} appeared before the present paper but nevertheless it should be considered a sequel to the present paper. In~\cite{agrcrdyn}, a variant of TCD maps on an infinite cylinder, called \emph{twisted TCD maps}, were introduced as a common generalization of twisted polygons~\cite{schwartz} and planar TCD maps. The connection to the bipartite dimer model and its cluster integrability~\cite{gkdimers} was then used to prove cluster integrability of cross-ratio dynamics as well as to provide an alternate proof in the case of the pentagram map.

The notion of VRC is dual to Lam's \emph{relation space}~\cite{LamCDM} (see also~\cite{abenda,abendagrinevich1,abendagrinevich2}). Our framework also has a similar flavor to that of \cite{fominpylyavskyy} where they introduce a master theorem in incidence geometry from which other theorems are derived. Each theorem is encoded by a bipartite graph in which every white vertex carries a point and every black vertex carries a line. However, their incidence conditions are different from ours and so far the relation between their work and ours is not understood.

\subsection*{Organization of the paper}

Section~\ref{sec:background} recalls the necessary background on projective geometry, triple crossing diagrams, bipartite graphs and vector-relation configurations. In Section~\ref{sec:tcdmapsvrc}, we define TCD maps and explain how they relate to VRCs. Section~\ref{sec:space} provides a constructive parametrization of TCD maps which is used to describe the moduli space of TCD maps. In Section~\ref{sec:tcddskp} we study local transformations of TCD maps and the connection to the dSKP equation. Section~\ref{sec:tcdconsistency} establishes multi-dimensional consistency and links the framework to classical integrable lattices. Section~\ref{sec:clusterstructures} introduces the projective and affine cluster structures. Section~\ref{sec:sections} introduces the operation of taking a section of a TCD map and establishes its basic properties. Finally, Section~\ref{sec:sectioncluster} proves the main correspondence theorem linking projective and affine cluster structures, completing the theoretical framework.

\section{Background}
\label{sec:background}

In this section, we recall the necessary background on projective geometry, triple crossing diagrams, black-trivalent bipartite graphs and vector-relation configurations.

\subsection{Projective and affine geometry}

To simplify exposition, we work over $\C$ in this paper. However, all results and definitions remain valid over $\R$ (and, in fact, over any field, though in some cases the statements may become trivial).

For $d\geq1$, \emph{projective space} $\CP^\dimension$ is defined as the space of all lines through the origin in $\C^{\dimension+1}$. A point in $\CP^\dimension$ is represented by an equivalence class
\[
	[x]=[x_1, \cdots, x_{\dimension+1}],
\]
where $x=(x_1, \dots, x_{\dimension+1}) \in \C^{\dimension+1} \setminus \{0\}$, and two vectors are identified if they are linearly dependent, that is, 
\[
	(x_1, \dots, x_{\dimension+1}) \sim \lambda (x_1, \dots, x_{\dimension+1}) \quad \text{for all } \lambda \in \C^* = \C \setminus\{0\}.
\]

A \emph{projective transformation} is an invertible linear transformation of $\C^{\dimension+1}$ modulo scalars, that is, an element of \[
	\mathrm{PGL}(\dimension+1,\C) = \mathrm{GL}(\dimension+1,\C)/\C^*. 
\]

An \emph{affine chart} on $\CP^\dimension$ is given by choosing a hyperplane $\hyperplane \subset \CP^\dimension$ and identifying the complement $\CP^\dimension \setminus \hyperplane$ with affine space $\C^{\dimension}$ via a choice of coordinates. We assume in the following that $\hyperplane$ is the hyperplane $x_{\dimension+1}=0$, giving the \emph{standard affine chart}
\[
	\{[x] \in \CP^\dimension: x_{\dimension+1} \neq 0\}
\]
which is identified with $\C^{\dimension}$ via the map
\[
	[x_1, \cdots, x_{\dimension+1}] \longmapsto \left( \frac{x_1}{x_{\dimension+1}}, \dots, \frac{x_{\dimension}}{x_{\dimension+1}} \right).
\]

An \emph{affine transformation} is an element of $\mathrm{PGL}(\dimension+1,\C)$ of the form 
\[
\begin{bmatrix}
A & b \\
0 & 1
\end{bmatrix} \quad \text{mod } \C^*,
\]
with $A\in\mathrm{GL}(\dimension,\C)$ and $b$ a $d$-dimensional column vector. The subgroup of affine transformations is precisely the stabilizer of the hyperplane $\hyperplane$. In other words, affine transformations preserve the standard affine chart.

Let $\centerofprojection\in \CP^\dimension$ be a point called the \emph{center} and let $\hyperplane \subset \CP^\dimension$ be a hyperplane not containing $\centerofprojection$. The \emph{central projection} from $\centerofprojection$ onto $\hyperplane$ is the rational map
\[
	\pi_{\centerofprojection,\hyperplane} : \CP^\dimension \dashrightarrow \hyperplane,
\]
defined by sending a point $\point \in \CP^\dimension \setminus \{\centerofprojection\}$ to the intersection point with $\hyperplane$ of the line joining $\centerofprojection$ to $\point$. More generally, we define central projection from a subspace as follows.

\begin{definition}[Central projection]\label{def:projection}
    Let $\subspacea\subset\CP^{\dimension}$ be a projective subspace of dimension $\dimension-r-1$ (the \emph{center}),
    and let $\subspaceb\subset\CP^{\dimension}$ be a complementary subspace of dimension $r$ (i.e. $\subspacea\cap\subspaceb=\varnothing$). The \emph{central projection} $\pi_{\subspacea}$ from $\subspacea$ onto $\subspaceb$ is the rational map
    \[
      \pi_{\subspacea}:\ \CP^\dimension \dashrightarrow \subspaceb \cong \CP^r,
    \]
    that sends every point $\point\in\CP^\dimension\setminus\subspacea$ to the intersection of $\subspaceb$ with the projective span of $\subspacea$ and $\point$.    
\end{definition}

Since we will quotient by the action of $\operatorname{PGL}(r{+}1)$ on the target, the specific choice of $\subspaceb$
and of the identification $\subspaceb\cong\CP^r$ is irrelevant; hence we write $\pi_{\subspacea}$ instead of $\pi_{\subspacea,\subspaceb}$.
It is also useful to describe central projections in homogeneous coordinates. Specifically, if $\tilde{\subspacea}\subset\C^{\dimension+1}$ is the linear subspace corresponding to $\subspacea$, then $\pi_{\subspacea}$ is the projectivization of the linear quotient
\begin{equation}\label{eq:central_projection_affine}
\tilde \pi_{\tilde \subspacea}: \C^{\dimension+1} \longrightarrow \C^{\dimension+1}/\tilde{\subspacea} \cong \C^{r+1}.
\end{equation}

The following geometric ratios will be used to define cluster variables associated to TCD maps.

\begin{definition}[Oriented length ratio]\label{def:orientedlengthratio}
	Let $\point_1,\point_2,\point_3$ be points on a line $\projectiveline$ in $\CP^\dimension$ and assume we have a fixed affine chart $\C^{\dimension}$ of $\CP^\dimension$ that contains $\point_1,\point_2,\point_3$. The \emph{oriented length ratio} $\lambda(\point_1,\point_2,\point_3)$ is defined as 
	\begin{align*}
		\lambda(\point_1,\point_2,\point_3) = \frac{\point_1-\point_3}{\point_2-\point_3},
	\end{align*}
	where subtraction and division are in the affine chart $\C \subset \projectiveline$ induced by the intersection $\projectiveline\cap\C^\dimension$. 
\end{definition}

The oriented length ratio is an affine invariant. However, it depends on the choice of affine chart and is not invariant under general projective transformations.

\begin{definition}[Multi-ratio]\label{def:multiratio}
	Let $\point_1,\point_{1,2},\point_2,\point_{2,3},\dots, \point_m,\point_{m,1}$ be $2m$ points in $\CP^\dimension$ such that each point $\point_{i,i+1}$ is on the line ${\point_i \point_{i+1}}$ (indices taken modulo $m$). Choose an affine chart $\C^{\dimension} \subset \CP^\dimension$ containing all the points and define the \emph{multi-ratio} of these points as
	\begin{align*}
		\mr(\point_1,\point_{1,2},\point_2,\dots, \point_m,\point_{m,1}) = \prod_{i=1}^m \frac{\point_i - \point_{i,i+1}}{\point_{i,i+1} - \point_{i+1}} = (-1)^m \prod_{i=1}^m \lambda(\point_i,\point_{i+1},\point_{i,i+1}).
	\end{align*}

	In the special case $m=2$, the multi-ratio becomes the classical \emph{cross-ratio}:
\[
		\cro(\point_1,\point_{1,2},\point_2,\point_{2,1}) = \frac{\point_1 - \point_{1,2}}{\point_{1,2} - \point_{2}}\frac{\point_2 - \point_{2,1}}{\point_{2,1} - \point_{1}}. \qedhere
\]
\end{definition}

Let us gather a few basic but important properties of the multi-ratio; see, for example, \cite{ddgbook} for proofs.
\begin{lemma}\label{lem:projinvariants}
	Let $\point_1,\point_{1,2},\point_2,\point_{2,3},\dots, \point_m,\point_{m,1}$ be $2m$ points in $\CP^\dimension$ such that every $\point_{i,i+1}$ lies on the line $\point_i \point_{i+1}$ (indices taken modulo $m$). Then the multi-ratio $\mr(\point_1,\point_{1,2},\point_2,\dots, \point_m,\point_{m,1})$ satisfies the following properties:
	\begin{enumerate}
		\item It is independent of the choice of affine chart used to compute it.
		\item It is invariant under projective transformations.
		\item It is invariant under central projections, provided the center does not meet $\point_i \point_{i+1}$ for all $i$. \qedhere
	\end{enumerate}
\end{lemma}

\subsection{Triple crossing diagrams}\label{sec:tcds}

   What we refer to as a \emph{triple crossing diagram} was originally introduced by Dylan Thurston \cite{thurstontriple} under the names \emph{triple-point diagram} and \emph{triple diagram}. These diagrams are also sometimes called \emph{Thurston diagrams} \cite{fmloop}, and the term {triple crossing diagram} has appeared in the literature in \cite{bocklandtdimerabc,bwtriple}. In what follows, we recall several basic definitions and results from Thurston’s paper, before relating triple crossing diagrams to vector-relation configurations in the next section.

\begin{definition}[Triple crossing diagrams]
	A \emph{triple crossing diagram} or \emph{TCD} for short is a collection $\tcd$ of oriented closed intervals and oriented circles immersed smoothly into a disk $\disk$. The image of a connected component is called a \emph{strand}; it is either an \emph{arc} (the image of an interval) or a \emph{loop} (the image of a circle). The immersions are required to satisfy:
	\begin{enumerate}
		\item Exactly three strands cross at each point of intersection.
		\item The endpoints of arcs are distinct points on the boundary of the disk, and no other points map to the boundary.
		\item The orientations of the strands induce consistent orientations on faces of $\tcd$ (that is, on the complementary regions of $\tcd$).
	\end{enumerate}
  Let $\tcdm$ (resp.~$\tcdp$) denote the set of clockwise (resp.~counterclockwise) faces of $\tcd$, including the boundary faces. 
	\end{definition}

    TCDs are considered up to ambient isotopies of $\disk$. Due to Condition~(3), the endpoints of strands have to alternate in orientation between sources and sinks around the boundary of the disk. Suppose there are $n$ strands that are arcs, so there are $n$ endpoints that are sources and $n$ that are sinks. We label the sources (resp.~sinks) by ${1},\dots,{n}$ (resp.~$\overline{1},\dots, \overline{n}$) in clockwise order such that source ${i}$ is between sinks $\overline{i}$ and $\overline{i+1}$. We label arcs by the labels of their sources.

\begin{definition}[Strand permutation of a TCD]\label{def:endpointmatching}
	Let $\tcd$ be a TCD with $n$ strands that are arcs. Let $C_\tcd \in S_n$ be the permutation defined by \[C_\tcd(i) = j \text{ if there is a strand from source $i$ to sink $\overline j$}.\] 
    We call $C_\tcd$ the \emph{strand permutation} of $\tcd$.
    If $1\leq k\leq n$, we denote by $\enm nk$ the permutation sending any $i$ to $i+k \mod n$.
\end{definition}

\begin{definition}[Minimal TCDs]\label{def:minimaltcd}
    A \emph{minimal} TCD is a TCD with no loops and the fewest triple crossings among all TCDs with the same strand permutation.
\end{definition}

Minimal TCDs can be characterized by the absence of certain forbidden configurations of strands.

\begin{theorem}[\protect{\cite[Theorem~7]{thurstontriple}}]\label{th:mintcdforbidden}
	A TCD is minimal if and only if
	\begin{enumerate}
		\item No strand is a loop.
		\item No strand has a self-intersection.
		\item No pair of strands forms a parallel bigon, that is, there are no two strands $\alpha,\beta$ and two triple crossings $x$~and~$y$ with both strands oriented from $x$ to $y$.\qedhere
	\end{enumerate}
\end{theorem}

For the remainder of this paper, we will assume without further mention that \emph{all TCDs are minimal}.

\begin{theorem}[\protect{\cite[Theorem~3]{thurstontriple}}]
    All $n!$ strand permutations are realizable by minimal TCDs.
\end{theorem}

Explicit constructions of minimal TCDs realizing any given strand permutation appear in~\cite{thurstontriple}~and~\cite{postgrass}.

\begin{figure}
	\centering
	\raisebox{-0.5\height}{\includegraphics[scale=0.7,angle=90]{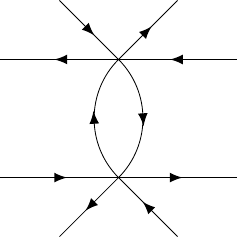}}
	\hspace{2mm}$\longleftrightarrow$\hspace{2mm}
	\raisebox{-0.5\height}{\includegraphics[scale=0.7]{tikz/tcdeye}}
	\hspace{15mm}
	\raisebox{-0.5\height}{\includegraphics[scale=0.7,angle=90]{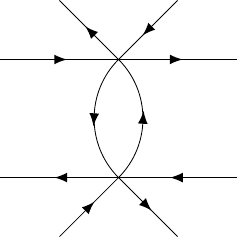}}
	\hspace{2mm}$\longleftrightarrow$\hspace{2mm}
	\raisebox{-0.5\height}{\includegraphics[scale=0.7]{tikz/tcdeyemirror}}
	\caption{The 2-2 move: for a clockwise bigon (left) and a counterclockwise bigon (right).}\label{fig:twotwo}
\end{figure}

\begin{definition}[2-2 move]
	A \emph{2-2 move} $s:\tcd \rightsquigarrow \tilde\tcd$ is one of the two local moves shown in Figure~\ref{fig:twotwo}, replacing a portion of the TCD $\tcd$ with the same portion rotated by an angle of $\pi/2$. We say that two minimal TCDs are \emph{move-equivalent} if they are related by a sequence of 2-2 moves. 
\end{definition}

The 2-2 move plays a central role in the geometric dynamics that we study. Note that the orientation of the bigon is preserved under the move. Since we distinguish between orientations, we will treat the 2-2 move applied to a clockwise bigon as different from the move applied to a counterclockwise bigon. One will correspond to a change in local geometry, while the other will be a reparametrization (see Section~\ref{sec:tcddskp}).

Finally, we have the following simple characterization of move-equivalence.

\begin{theorem}[\protect{\cite[Theorem~5]{thurstontriple}}]\label{th:tcdmovesconnected}
	Any two minimal TCDs with the same strand permutation are move-equivalent.
\end{theorem}

\subsection{Black-trivalent bipartite graphs}

The following notion was introduced by Lam~\cite{Lam} in the closely related context of electrical networks. On a first reading, the reader may assume that all graphs lie in a disk~$\disk$; the more technical notion of a cactus is only needed to handle degenerate examples, that will regularly arise in Section~\ref{sec:sections} when discussing some transformations called ``sections''.

\begin{definition}[Cacti]
    Consider a disk $\disk$ with $n$ white vertices on its boundary labeled $w_1^\partial,\dots,w_n^\partial$ in clockwise cyclic order. A \emph{noncrossing partition} $\ncp$ of $\{w_1^\partial,\dots,w_n^\partial\}$ is a partition of $\{w_1^\partial,\dots,w_n^\partial\}$ into disjoint nonempty blocks such that there are no two distinct blocks $B_1$ and $B_2$ with $w_{i_1}^\partial, w_{i_3}^\partial \in B_1$ and $w_{i_2}^\partial,w_{i_4}^\partial \in B_2$ and $i_1<i_2<i_3<i_4$. 

    A \emph{cactus} is a topological space $\disk/\ncp$ obtained by gluing the boundary white vertices according to a non-crossing partition $\ncp$ of $\{w_1^\partial,\dots,w_n^\partial\}$. It is a union of disks glued together along some white vertices on their boundaries.
\end{definition}

\begin{figure}
	\centering
	\raisebox{-0.5\height}{\includegraphics[scale=1]{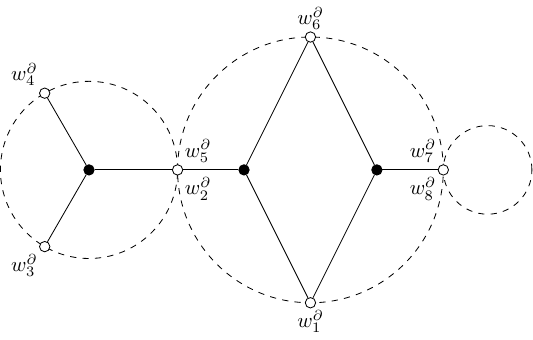}}	
	\caption{A black-trivalent bipartite graph $\pb$.}
	\label{fig:\plabic{}_graph}
\end{figure}

We now define a certain class of graphs that will be equivalent to TCDs. The triple crossing condition is reflected in a degree condition on black vertices.

\begin{definition}[Black-trivalent bipartite graphs]
    A \emph{black-trivalent bipartite graph (BTB graph for short)} $\pb = (B \sqcup W,E,F)$ is a planar bipartite graph embedded in a cactus with boundary vertices $w_1^\partial, \dots, w_n^\partial$ such that every black vertex is trivalent. 
\end{definition}

We fix the following notation for \plabic{} graphs. We let $W^\partial$ and $F^\partial$ denote the sets of boundary white vertices and boundary faces respectively. The set of 
internal white vertices and faces we denote by
\begin{align*}
    W^\circ:= W \setminus W^\partial, \quad \mbox{and} \quad F^\circ:=F\setminus F^\partial,
\end{align*}
respectively. Moreover, we use the convention that $f_{i,i+1}^\partial$ is the boundary face between $w_i^\partial$ and $w_{i+1}^\partial$.

\Plabic{} graphs are a modification of Postnikov's plabic graphs \cite{postgrass} chosen so as to be dual to triple crossing diagrams. Other similar versions have appeared in the literature under the names white-partite~\cite{GPW} and black-trivalent~\cite{galashincritical}. 

\begin{example}
	Figure~\ref{fig:\plabic{}_graph} shows a \plabic{} graph in a cactus $\disk/\pi$, where \[\pi=\{1\}\sqcup\{2,5\}\sqcup\{3\}\sqcup\{4\}\sqcup\{6\}\sqcup\{7,8\}. \qedhere\] 
\end{example}

\begin{figure}
	\centering
	\begin{tabular}{c@{\qquad}c@{\qquad}c}	

		\raisebox{-0.5\height}{\includegraphics[scale=0.72]{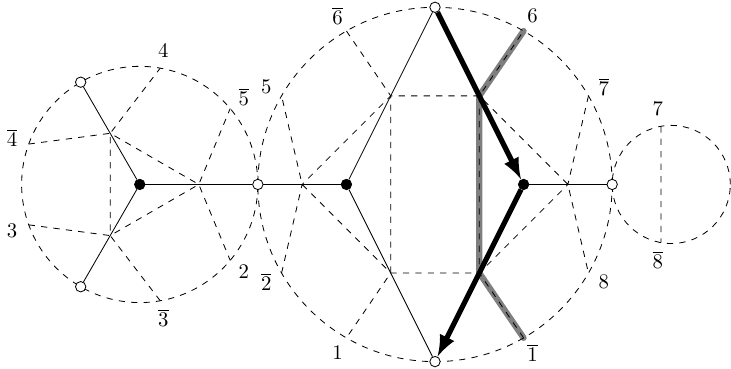}}			
		&
		\(\longrightarrow\)
		&		

		\raisebox{-0.5\height}{\includegraphics[scale=0.65]{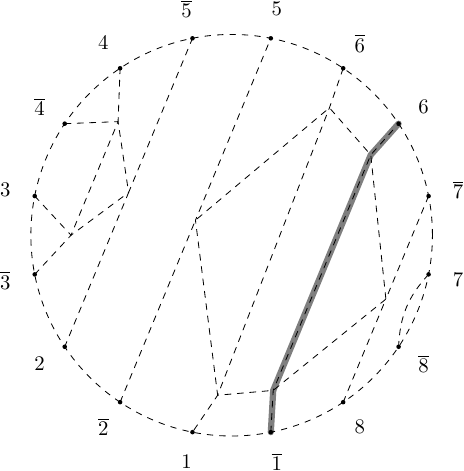}}		
	\end{tabular}	
	\caption{The construction of the medial graph $\pb^\times$ of the graph $\pb$ in Figure~\ref{fig:\plabic{}_graph}. A medial strand (bold gray) and the corresponding zigzag path (bold black) are highlighted.}
	\label{fig:medial_graph}
\end{figure}

\begin{definition}[Medial graph $\pb^\times$ of $\pb$]
    Let \( \pb^\times \) denote the \emph{medial graph} of \( \pb \) obtained from \( \pb \) as follows. Place vertices $\overline{1},1,\dots,  \overline{n}, n$ in clockwise order around the boundary of the cactus such that $w_i^\partial$ is between $\overline{i}$ and $i$ and a vertex $v^\times_e$ in the middle of each edge $e$ of $\pb$. Connect $v^\times_e$ to $v^\times_{e'}$ if they occur consecutively around the boundary of a face of $\pb$. Draw an edge from $i$ (resp.~$\overline{i}$) to $v_e^\times$ if $e$ is the last (resp.~first) edge in clockwise order around the boundary face of $\pb$ containing $i$ (resp.~$\overline i$). If $w_i^\partial$ is an isolated vertex, draw an edge between $i$ and $\overline{i}$. By construction, each internal vertex of \(G^\times\) has degree \(4\). Finally take the preimage of this graph under the quotient map from $\disk$ to $\disk/\ncp$ to get $\pb^\times$. 
\end{definition}

Note that, while the graph $G$ lives on a cactus, its medial graph $G^\times$ lives on the disk. The faces of $\pb^\times$ are of three types, corresponding to black vertices, white vertices, and faces of $\pb$. 

\begin{example}
    The medial graph of the graph $\pb$ in Figure~\ref{fig:\plabic{}_graph} is shown on the right picture of Figure~\ref{fig:medial_graph}.
\end{example}

\begin{definition}[Medial strands, zig-zag paths and minimality] \label{def:zigzags}
    A \emph{medial strand} of $\pb^\times$ is a path in $\pb^\times$ that starts at a boundary vertex $i$ and goes straight through at each (degree 4) internal vertex, i.e. the path has one edge on its left and one edge on its right at each internal vertex.
    
    The \emph{zig-zag path} of $\pb$ corresponding to a medial strand is the path in $\pb$ that uses edges in the same sequence as the medial strand; such a path turns maximally left at white vertices and maximally right at black vertices. At boundary white vertices, a zig-zag path always begins along the leftmost edge and terminates when it exits the graph through the rightmost edge.
    
     We say that $\pb$ is \emph{minimal} or \emph{reduced} if $G^\times$ has no medial strands that form closed loops, no self-intersecting medial strands, and no pair of medial strands that forms a \emph{parallel bigon}, that is, there are no two medial strands $\alpha,\beta$ that go through two distinct vertices $v_{e_1}^\times,v_{e_2}^\times$ with both of them oriented from $v_{e_1}^\times$ to $v_{e_2}^\times$. 
\end{definition}

A consequence of Definition~\ref{def:zigzags} is that a minimal \plabic{} graph has no loops or multiple edges, and no isolated internal vertices; see \cite{postgrass}. Hereafter, we assume that all our \plabic{} graphs are minimal.

\begin{definition}[Strand permutation of a \plabic{} graph]
    The \emph{strand permutation} $C_\pb \in S_n$ of a \plabic{} graph $G$ with $n$ boundary vertices is defined by
    \[
    C_\pb(i) = j \text{ if there is a medial strand from $i$ to $\overline{j}$.} \qedhere
    \]
\end{definition}

It is easy to see that $i$ is a fixed point of $C_\pb$ if and only if $w_i^\partial$ is an isolated vertex and forms a singleton block of $\ncp$.

\begin{figure}
    \centering

	\raisebox{-0.5\height}{\includegraphics[scale=1]{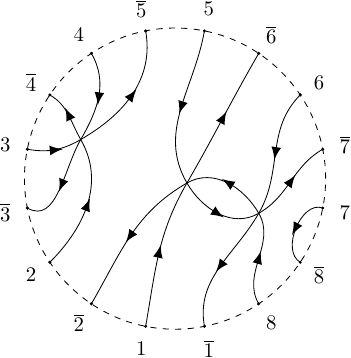}}
    \caption{The TCD of the \plabic{} graph $\pb$ from Figure~\ref{fig:\plabic{}_graph}.
    }
    \label{fig:tcd_from_\plabic{}_graph}
\end{figure}

The following definition explains the correspondence between BTB graphs and TCDs.

\begin{definition}[The \plabic{} graph of a TCD]\label{def:pb}
     Let $\tcd$ be a TCD in a disk $\disk$.  The associated \plabic{} graph $\pb$ is constructed as follows:
    \begin{itemize}
        \item Place a white vertex $w$ in each interior clockwise face of $\tcd$. 
        \item For each $i$, place a boundary white vertex $w_i^\partial$ between sink $\overline{i}$ and source $i$. 
        \item For each clockwise boundary face, glue all the boundary white vertices $w_i^\partial$ that lie in that face into a single white vertex, possibly creating a cactus.
        \item Place a black vertex $b$ at each triple crossing of $\tcd$.
        \item Draw an edge between $b$ and $w$ if the triple crossing at $b$ lies on the boundary of the face containing $w$.
    \end{itemize}
    Conversely, we obtain the TCD $\tcd$ from the \plabic{} graph $\pb$ by first drawing the medial graph $\pb^\times$, orienting medial edges clockwise around white vertices and counterclockwise around black vertices, and then contracting the triangular faces corresponding to black vertices of $\pb$.
\end{definition}
This construction yields a natural correspondence between TCDs~$\tcd$ and \plabic{} graphs~$\pb$. Under this correspondence,  

\begin{align*}
    \mbox{clockwise faces of }\tcd  &\longleftrightarrow   \mbox{white vertices of }\pb, \\
    \mbox{triple crossings of }\tcd  &\longleftrightarrow \mbox{black vertices of }\pb, \\
    \mbox{counterclockwise faces of }\tcd  &\longleftrightarrow  \mbox{faces of }\pb.
\end{align*}

Moreover, strands of $\tcd$ correspond to medial strands of $\pb^\times$, the notions of minimality coincide, and $C_\tcd = C_\pb$.

\begin{example}
    Figure~\ref{fig:tcd_from_\plabic{}_graph} shows the TCD corresponding to the \plabic{} graph from Figure~\ref{fig:\plabic{}_graph}.
\end{example}

\subsection{Vector-relation configurations}

In this subsection, we recall the notion of vector-relation configurations following~\cite{vrc}. Recall that a set of vectors is called a \emph{circuit} if they are a minimal linearly dependent set, that is, if they are linearly dependent but any proper subset is linearly independent.

\begin{definition}[Vector-relation configurations]\label{def:vrc}
	Let $d\geq1$ be an integer and let $\pb = (B \sqcup W, E, F)$ be a planar bipartite graph embedded in a cactus. A \emph{vector-relation configuration} on $\pb$ (or \emph{VRC} for short) is a pair $(\vrc, \relation)$, where:
\begin{enumerate}
    \item $\vrc: W \to \C^{\dimension+1}$ is an assignment of a nonzero vector $\vrc(w)$ to each white vertex $w \in W$.
   
    \item For each black vertex $b \in B$, the set of vectors $\{\vrc(w) : bw \in E\}$ attached to white neighbors of $b$ satisfy the linear dependence relation
    \[
    \relation(b):\quad \sum_{bw \in E} \edgeweight(bw)\, \vrc(w) = 0, 
    \]
    where the coefficients $\edgeweight(bw) \in \C$ are called \emph{edge weights}. 
    
    We call a VRC $(\vrc, \relation)$ a \emph{circuit configuration} if, in addition, for each black vertex $b \in B$, the set of vectors 
\(
\{\vrc(w) : bw \in E\}
\) 
forms a circuit. The circuit condition implies that $\vrc(w)$ and $\edgeweight(bw)$ never vanish. \qedhere
\end{enumerate}
\end{definition}

In the definition of VRCs, we do not assume that black vertices are trivalent. Accordingly, we use the term planar bipartite graph for general graphs of this type.

\begin{definition}[Gauge transformations for VRCs] 
    A \emph{gauge transformation} of a VRC $(\vrc, \relation)$ is a rescaling of its vectors, relations and edge weights determined by a function $\lambda : B \sqcup W \longrightarrow \C^*$ as follows:
\begin{itemize}
    \item Each vector is rescaled by $\vrc(w) \longmapsto \lambda(w)\,\vrc(w)$.
    \item Each relation is rescaled by $\relation(b) \longmapsto \lambda(b)\,\relation(b)$.
    \item Each edge weight transforms as $
    \edgeweight(bw) \longmapsto \lambda(b) \lambda(w)^{-1} \edgeweight(bw)$. \qedhere
\end{itemize}
\end{definition}

Such gauge transformations are classical in the study of the dimer model in statistical mechanics, see for example \cite{kenyondimerintro}.

\newlength{\FigColW}\setlength{\FigColW}{0.38\textwidth} 
\newlength{\FigMidW}\setlength{\FigMidW}{2.6em}          
\newlength{\FigLblW}\setlength{\FigLblW}{1.2cm}          
\newcommand{\AlignedRow}[2]{
  \begin{tabular}{@{}c c c@{}}
    \makebox[\FigColW][c]{
      \begin{tikzpicture}[baseline={(current bounding box.center)}]
        #1
      \end{tikzpicture}
    } &
   
    \makebox[\FigMidW][c]{
      \begin{tikzpicture}[baseline={(current bounding box.center)}]
        \node {$\longleftrightarrow$};
      \end{tikzpicture}
    } &
 
    \makebox[\FigColW][c]{
      \begin{tikzpicture}[baseline={(current bounding box.center)}]
        #2
      \end{tikzpicture}
    }
  \end{tabular}
}

\newcommand{\RowWithLabel}[3]{
  \begin{tabular}{@{}c c@{}}
    \makebox[\FigLblW][c]{#1} &
    \makebox[0.95\textwidth][c]{\AlignedRow{#2}{#3}}
  \end{tabular}
}

\begin{figure}
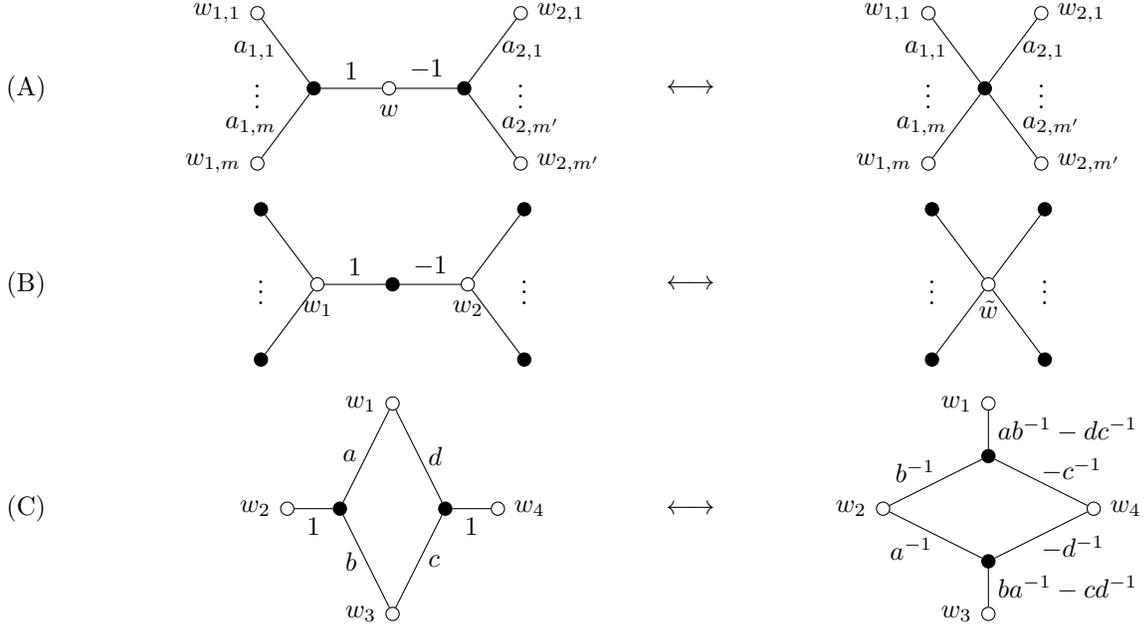

\centering

\begin{subfigure}{\textwidth}
  \centering
  \RowWithLabel{(A)}{
    \coordinate (fig) at (0,0);
    \node[wvert, label=below:$w$] (p) at (0,0) {};
    \node[bvert] (b1) at (-1,0) {};
    \node[bvert] (b2) at (1,0) {};

    \coordinate[wvert,label=left:$w_{1,1}$] (l1) at (-1.75,1) {};
    \node at (-1.75,0) {$\vdots$};
    \coordinate[wvert,label=left:$w_{1,m}$] (l4) at (-1.75,-1) {};

    \coordinate[wvert,label=right:$w_{2,1}$] (r1) at (1.75,1) {};
    \node at (1.75,0) {$\vdots$};
    \coordinate[wvert,label=right:$w_{2,m'}$] (r4) at (1.75,-1) {};

    \draw (l1) -- node[left] {$a_{1,1}$} (b1);
    \draw (b1) -- node[above] {$1$} (p);
    \draw (p)  -- node[above] {$-1$} (b2);
    \draw (b1) edge node[left]  {$a_{1,m}$} (l4);
    \draw (b2) edge node[right] {$a_{2,1}$} (r1);
    \draw (b2) edge node[right] {$a_{2,m'}$} (r4);
  }{
    \coordinate (fig) at (0,0);
    \node[bvert] (b) at (0,0) {};

    \coordinate[wvert,label=left:$w_{1,1}$] (l1) at (-.75,1) {};
    \node at (-.75,0) {$\vdots$};
    \coordinate[wvert,label=left:$w_{1,m}$] (l4) at (-.75,-1) {};

    \coordinate[wvert,label=right:$w_{2,1}$] (r1) at (.75,1) {};
    \node at (.75,0) {$\vdots$};
    \coordinate[wvert,label=right:$w_{2,m'}$] (r4) at (.75,-1) {};

    \draw (l1) -- node[left] {$a_{1,1}$} (b);
    \draw (b)  edge node[left]  {$a_{1,m}$} (l4);
    \draw (b)  edge node[right] {$a_{2,1}$} (r1);
    \draw (b)  edge node[right] {$a_{2,m'}$} (r4);
  }
\end{subfigure}

\vspace{2mm}

\begin{subfigure}{\textwidth}
  \centering
  \RowWithLabel{(B)}{

    \coordinate (fig) at (0,0);
    \node[bvert] (p) at (0,0) {};
    \coordinate[wvert,label=below:$w_1$] (wL) at (-1,0) {};
    \coordinate[wvert,label=below:$w_2$] (wR) at (1,0) {};

    \coordinate[bvert] (l1) at (-1.75,1) {};
    \node at (-1.75,0) {$\vdots$};
    \coordinate[bvert] (l4) at (-1.75,-1) {};

    \coordinate[bvert] (r1) at (1.75,1) {};
    \node at (1.75,0) {$\vdots$};
    \coordinate[bvert] (r4) at (1.75,-1) {};

    \draw (l1) -- (wL);
    \draw (wL) -- node[above] {$1$} (p);
    \draw (p)  -- node[above] {$-1$} (wR);
    \draw (wL) edge (l4);
    \draw (wR) edge (r1);
    \draw (wR) edge (r4);
  }{
    \coordinate (fig) at (0,0);
    \node[wvert,label=below:$\tilde w$] (w) at (0,0) {};

    \coordinate[bvert] (l1) at (-.75,1) {};
    \node at (-.75,0) {$\vdots$};
    \coordinate[bvert] (l4) at (-.75,-1) {};

    \coordinate[bvert] (r1) at (.75,1) {};
    \node at (.75,0) {$\vdots$};
    \coordinate[bvert] (r4) at (.75,-1) {};

    \draw (l1) -- (w);
    \draw (w) edge (l4);
    \draw (w) edge (r1);
    \draw (w) edge (r4);
  }
\end{subfigure}

\vspace{2mm}

\begin{subfigure}{\textwidth}
  \centering
  \RowWithLabel{(C)}{
    \begin{scope}[scale=.7]
      \coordinate (fig) at (0,0);
      \node[wvert,label=left:$w_2$] (w2) at (0,0) {};
      \node[bvert] (b2)  at (1,0) {};
      \node[wvert,label=left:$w_1$] (w1)  at (2,2) {};
      \node[wvert,label=left:$w_3$] (w3) at (2,-2) {};
      \node[bvert] (b4) at (3,0) {};
      \node[wvert,label=right:$w_4$] (w4) at (4,0) {};
      \path[-]
        (w2) edge node[below] {$1$} (b2)
        (b2) edge node[left] {$a$} (w1)
        (b2) edge node[left] {$b$} (w3)
        (w3) edge node[right] {$c$} (b4)
        (w1) edge node[right] {$d$} (b4)
        (b4) edge node[below] {$1$} (w4);
    \end{scope}
  }{
    \begin{scope}[scale=.7]
      \coordinate (fig) at (0,0);
      \node[wvert,label=left:$w_2$] (w2) at (0,0) {};
      \node[bvert] (b1)  at (2,1) {};
      \node[wvert,label=left:$w_1$] (w1)  at (2,2) {};
      \node[wvert,label=left:$w_3$] (w3) at (2,-2) {};
      \node[bvert] (b3) at (2,-1) {};
      \node[wvert,label=right:$w_4$] (w4) at (4,0) {};
      \path[-]
        (b1) edge node[right] {$ab^{-1}-dc^{-1}$} (w1)
        (b1) edge node[above left,inner sep=0] {$b^{-1}$} (w2)
        (b1) edge node[above right,inner sep=0] {$-c^{-1}$} (w4)
        (b3) edge node[right] {$ba^{-1}-cd^{-1}$} (w3)
        (b3) edge node[below left,inner sep=1] {$a^{-1}$} (w2)
        (b3) edge node[below right,inner sep=0] {$-d^{-1}$} (w4);
    \end{scope}
  }
\end{subfigure}

\caption{Local moves for VRCs: (A) and (B) show contractions (left to right) and splits (right to left), respectively; (C) shows a spider move. The vertical dots indicate that several edges of the same type are attached.}
\label{fig:vrc_moves}
\end{figure}

\begin{definition}[Moves for VRCs]
A \emph{move} of a VRC relates a configuration $(\vrc, \relation)$ to another configuration $(\tilde{\vrc}, \tilde{\relation})$, obtained by locally replacing one subconfiguration with another.  
The three types of local moves for VRCs, shown in Figure~\ref{fig:vrc_moves}, are:
\begin{enumerate}
    \item[(A)] \emph{Contraction} of a degree-$2$ white vertex (left to right), or a \emph{split} of a black vertex (right to left).  
    For the split, the new vectors satisfy
    \begin{align*}
        \vrc(w) &= -a_{1,1}\,\tilde{\vrc}(w_{1,1}) - \cdots - a_{1,m}\,\tilde{\vrc}(w_{1,m}).
    \end{align*}

    \item[(B)] \emph{Contraction} of a degree-$2$ black vertex (left to right), or a \emph{split} of a white vertex (right to left).  
    In this case, the new vectors satisfy
    \begin{align*}
        \tilde{\vrc}(\tilde{w}) &= \vrc(w_1) = \vrc(w_2).
    \end{align*}

    \item[(C)] \emph{Spider move} at a face of degree four.
\end{enumerate}
In each case, the relations $\relation$ and $\tilde{\relation}$ are determined by the edge weights, which are also shown in Figure~\ref{fig:vrc_moves}.
\end{definition}

The reader familiar with local moves for the dimer model in statistical mechanics may be surprised by the observation that in the split, the edge weights of the two edges adjacent to the central degree-two vertex are opposite, rather than equal (left side of (A) and (B) in Figure~\ref{fig:vrc_moves}). The reason is that the edge weights of a VRC play the role of the product of the positive dimer edge weights with the Kasteleyn signs in the classical planar bipartite dimer theory, and the Kasteleyn signs of these two edges are opposite; see \cite{agrcrdyn} for more details.

\section{Triple crossing diagram maps and vector-relation configurations}
\label{sec:tcdmapsvrc}

In this section, we introduce TCD maps and describe their relationship with VRCs.

\begin{definition}[Triple crossing diagram maps]\label{def:tcdmap}
    Let $d\geq1$ and let $\tcd$ be a TCD and $\pb$ its \plabic{} graph. A \emph{triple crossing diagram map} or \emph{TCD map} is a map 
    \[
    \tcdmap: \tcdm \longrightarrow \CP^\dimension,
    \]
    that assigns a point in $\CP^\dimension$ to each clockwise face of $\tcd$ such that at every triple crossing, the three points associated to the three adjacent clockwise faces are distinct and lie on a line.
\end{definition}

Equivalently, identifying clockwise faces with white vertices in the \plabic{} graph associated to $\tcd$, a TCD map is a map
\[
    \tcdmap: W \longrightarrow \CP^\dimension,
\]
assigning a point in $\CP^\dimension$ to each white vertex, such that if a black vertex $b$ is incident to white vertices $w_1, w_2, w_3$, then the points $\tcdmap (w_1),\tcdmap (w_2), \tcdmap (w_3)$ are distinct and lie on a line $\projectiveline_b$. 

Let us also fix notation for the moduli space of TCD maps.

\begin{definition}[Rank and moduli space of TCD maps]
The \emph{rank} of a TCD map~$\tcdmap$, denoted~$\mrank(\tcdmap)$, is the dimension of the span of its image points.  
We denote by
\[
  \maps(\tcd,\dimension) 
  := \Bigl\{ \text{TCD maps $\tcdmap$ on $\tcd$ with $\mrank(\tcdmap)=\dimension$} \Bigr\} 
  \big/ \operatorname{PGL}(\dimension+1),
\]
the \emph{moduli space of rank-$\dimension$ TCD maps on $\tcd$}.
\end{definition}

\begin{proposition}
    There is a one-to-one correspondence between TCD maps and gauge-equivalence classes of circuit VRCs.
\end{proposition}
\proof{
    Given a circuit configuration $(\vrc,\relation)$ on a \plabic{} graph $\pb$, we obtain a TCD map by setting 
    \[
        \tcdmap(w) := [\vrc(w)] \in \CP^\dimension,
    \] 
    for all white vertices $w$, where $[\vrc(w)]$ is the projectivization of $\vrc(w)$. If $w,w'$ are adjacent to the same black vertex, then $\tcdmap(w)$ and $\tcdmap(w')$ are different points due to the circuit condition. Gauge-equivalent circuit configurations give rise to the same TCD map, hence, $\tcdmap$ is well-defined.

    Conversely, let $\tcdmap$ be a TCD map. For each white vertex $w \in W$, choose a lift $\vrc(w) \in \C^{\dimension+1} \setminus \{0\}$ of the point $\tcdmap(w) \in \CP^\dimension$.    
     For each black vertex $b$ adjacent to three white vertices $w_1, w_2, w_3$, the condition that $\tcdmap(w_1),\tcdmap(w_2),\tcdmap(w_3)$ lie on a line is equivalent to the existence of a linear dependence relation
    \[
        \relation(b):\quad \alpha_1 \vrc(w_1) + \alpha_2 \vrc(w_2) + \alpha_3 \vrc(w_3) = 0, \qquad \text{for some $\alpha_1,\alpha_2,\alpha_3 \in \C$},
    \]
    on lifts. The condition that the points are distinct implies that $\{\vrc(w_1),\vrc(w_2),\vrc(w_3)\}$ forms a circuit and in particular that $\alpha_1,\alpha_2,\alpha_3$ are nonzero.  
    Since different choices of lifts and different rescalings of the relations are gauge transformations, this defines a VRC on $\pb$ uniquely up to gauge. \qed

}

We regard the TCD maps as the fundamental theoretical objects instead of the VRCs, since we are primarily interested in objects and maps in projective geometry. Many of the interesting structures we study later on, such as cluster structures and partition functions, are defined modulo gauge transformations of VRCs. 
Considering TCD maps avoids having to quotient out by gauge transformations. Furthermore, local moves have a simpler geometric interpretation for TCD maps than for VRCs, see Section~\ref{sec:tcddskp}. The bipartite graphs carrying VRCs may look more general, because their black vertices have arbitrary degree. However, contracting all the degree-two black vertices and repeatedly splitting degree-two white vertices to reduce the degree of black vertices that have more than three neighbors, as in Figure~\ref{fig:vrc_moves}, one obtains a bipartite graph with trivalent black vertices carrying a VRC equivalent to the original one.

\begin{remark}
    In general, we attempt to minimize calculations and rely on incidence theorems and a flexible collection of lemmas. However, on occasions where we need to resort to explicit calculations, working with VRCs is the most practical approach.
\end{remark}

There is a particular gauge that is useful both for explicit calculations and for the affine cluster structures that we introduce in Section \ref{sec:affcluster}.

\begin{definition}[Affine gauge]\label{def:affinegauge}
	Consider a VRC $(\vrc,\relation)$ on a \plabic{} graph $\pb$. The edge weights $\edgeweight$ of the VRC are said to be in \emph{affine gauge} if 
	\begin{align}
		\sum_{ bw \in E} \edgeweight(bw) = 0 \label{eq:affinegauge}
	\end{align}
	 for every $b\in B$.
\end{definition}

Let us explain why we call this affine gauge.

\begin{definition}[Generic hyperplane with respect to a TCD map]
    Let $\tcdmap$ be a TCD map. We say that a hyperplane $\hyperplane$ is \emph{generic with respect to $\tcdmap$} if it does not contain any of the points of $\tcdmap$.
\end{definition}

Suppose $\tcdmap$ is a TCD map and $\hyperplane$ is a hyperplane that is generic with respect to $\tcdmap$. By a projective transformation, assume that $\CP^\dimension \setminus \hyperplane = \C^\dimension$ is the standard affine chart. Let $(\vrc,\relation)$ be a VRC of $\tcdmap$. Then the $(\dimension+1)$-th coordinate of each $\vrc(w)$ is nonzero. By performing gauge transformations at the white vertices, we can normalize all $\vrc(w)$ so that their $(\dimension+1)$-th coordinate equals $1$. In this gauge, the edge weights satisfy the relation~\eqref{eq:affinegauge}. Moreover, this shows that the VRC of every TCD map admits an affine gauge.

\section{The space of TCD maps}
\label{sec:space}

\begin{figure}
	\begin{tikzpicture}
		\coordinate (a1) at (240:2);
		\coordinate[label=left:$a_3$] (aa1) at (60:2);
		\coordinate (a2) at (0:2);
		\coordinate[label=below:$a_2$] (aa2) at (180:2);
		\coordinate (a3) at (120:2);
		\coordinate[label=right:$a_1$] (aa3) at (300:2);
		\draw[]
			(aa1) edge[-latex] (a1)
			(aa2) edge[-latex] (a2)
			(aa3) edge[-latex] (a3)
		;
		\node[bvert] (b) at (0,0) {};
		\node[wvert] (w1) at (330:1.25) {};
		\node[wvert] (w2) at (90:1.25) {};
		\node[wvert] (w3) at (210:1.25) {};
		\draw[-]
			(b) edge[rorient] (w1) edge[orient] (w2) edge[rorient] (w3)
		;
	\end{tikzpicture}
	\caption{The orientation $\mathcal O$ at a black vertex of $\pb$.}
	\label{fig:liorientatblack}
\end{figure}
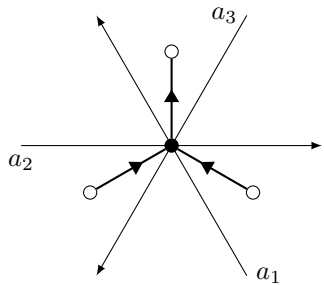

The goal of this section is to describe the moduli space of TCD maps. This relies on some combinatorial constructions that we now introduce.

\begin{definition}[Maximal \rank{}]
    Let $\tcd$ be a minimal TCD with \plabic{} graph $\pb$. The \emph{\rank{}} of $\tcd$ is defined as 
    \[
        \maxrank(\tcd) := |W| - |B| - 1. \qedhere
    \]
\end{definition}

The name will be justified in Proposition~\ref{prop:max_rank} where we show that $\maxrank(\tcd)$ is the maximal rank of any TCD map on $\tcd$. 

The following definition is due to Postnikov~\cite{postgrass}. A \emph{perfect orientation} on $\pb$ is an orientation such that each internal white vertex has a unique incoming edge and each black vertex has a unique outgoing edge.  
We define a particular perfect orientation $\mathcal O$ as follows; see Figure~\ref{fig:liorientatblack}.

\begin{definition}[The perfect orientation $\mathcal O$]
   Let $b\in B$ be a black vertex, and let $a_1,a_2,a_3$ be the three strands meeting at the triple crossing at $b$ with $a_1<a_2<a_3$. The cyclic order of strands around $b$ is the same as their cyclic order along the boundary (see, for example,~\cite[Lemma 3.4]{gkdimers}). For each edge $e=bw$ incident to $b$, in the \emph{perfect orientation} $\mathcal O$ we orient $e$ from $b$ to $w$ if $e$ is to the right of both $a_1$ and $a_3$; otherwise, orient $e$ from $w$ to $b$.
\end{definition}

Cyclically rotating the boundary labels yields $n$ distinct perfect orientations via the construction of $\mathcal O$. The choice of perfect orientation among these $n$ will not matter for us.

\begin{lemma}[Theorem~5.3 and Corollary~B.7 of \cite{mstwist}] \label{lem:perfectorientation}
The perfect orientation $\mathcal O$ satisfies the following properties:
\begin{enumerate}
    \item It is acyclic.
    \item It has $\maxrank(\tcd)+1$ sources, given by
    \[
      W^\partial \setminus \{ 1 \leq i \leq n  : C_{\tcd}(i) < i \}.\qedhere
    \]
\end{enumerate}
\end{lemma}

Since the perfect orientation $\mathcal O$ is acyclic, it defines a partial order on white vertices of $\pb$, which we will use to construct TCD maps inductively, vertex by vertex. 

\begin{definition}[$\mathcal O$-poset]
  We define the \emph{$\mathcal O$-poset} $\mathcal{P}_{\mathcal O} = (W, \le_{\mathcal O})$ to be the poset on $W$ induced by the orientation~$\mathcal O$.
\end{definition}

By induced, we mean that $w_1 \le_{\mathcal O} w_2$ if and only if there is an oriented path from $w_1$ to $w_2$ in $\mathcal O$.
By Lemma~\ref{lem:perfectorientation}, $\mathcal P_{\mathcal O}$ has $\maxrank(\tcd)+1$ minimal elements, given by
\[
    W^{\min}:=    W^\partial \setminus \{ 1 \leq i \leq n  : C_{\tcd}(i) < i \}.
\]

\begin{proposition}\label{prop:max_rank}
    Let $\tcd$ be a TCD. Every TCD map $\tcdmap$ on $\tcd$ satisfies
    \[
        \mrank(\tcdmap) \leq \maxrank(\tcd). \qedhere
    \]
\end{proposition}
\begin{proof}
    In the orientation $\mathcal O$, each non-minimal white vertex $w_1$ has a unique incoming edge $b \to w_1$. Let $w_2, w_3$ be the other two neighbors of $b$. Then $\tcdmap(w_1) \in \tcdmap(w_2) \tcdmap(w_3)$. Inductively, all white vertices lie in the span of the $\maxrank(\tcd)+1$ minimal ones, so the rank is at most $\maxrank(\tcd)$.
\end{proof}

Motivated by Proposition~\ref{prop:max_rank}, we say that a TCD map~$\tcdmap$ has \emph{maximal rank} if 
\[
  \mrank(\tcdmap) = \maxrank(\tcd).
\]
We now define an algorithm to construct TCD maps.

\begin{definition}[Construction algorithm]\label{def:construction_algorithm}
    Let $\tcd$ be a TCD with \plabic{} graph $\pb$. Let $1 \leq \dimension \leq \maxrank(\tcd)$. The following \emph{construction algorithm} yields a parameterization of $\maps(\tcd,\dimension)$.
    \begin{enumerate}
        \item Choose a linear extension $\varepsilon$ of the poset $\mathcal{P}_{\mathcal O}$.
        \item Assign points $\tcdmap(w_i^\partial)\in\CP^\dimension$ to the minimal elements $w_i^\partial$ of $\mathcal{P}_{\mathcal O}$ subject to the condition that these points span $\CP^\dimension$. 
        \item Let $w_1$ be the $\varepsilon$-smallest white vertex for which $\tcdmap(w_1)$ has not yet been assigned. Let $bw_1$ be the unique incoming edge at $w_1$. Let $w_2,w_3$ be the other two white neighbors of $b$, and let $f$ be the face incident to $b$ between $w_2$ and $w_3$; by minimality, $w_2$ and $w_3$ are distinct. Proceed as follows:
        \begin{enumerate}
            \item If $\tcdmap(w_2)=\tcdmap(w_3)$, then the algorithm fails. 
            \item If $\tcdmap(w_2) \neq \tcdmap(w_3)$, place $\tcdmap(w_1)$ on the line $\tcdmap(w_2)\tcdmap(w_3)$, chosen distinct from both $\tcdmap(w_2)$ and $\tcdmap(w_3)$.
        \end{enumerate}
        \item Repeat Step~(3) until all points have been assigned. \qedhere
    \end{enumerate}
\end{definition}

Every TCD map can be obtained via the construction algorithm, since the algorithm never fails for a given TCD map; in particular, in does not depend on the choice of linear extension $\varepsilon$.

\begin{proposition}
    The moduli space $\maps(\tcd,\dimension)$ is a variety of dimension 
    \[
        (\dimension-1)(\maxrank(\tcd)+1)-(\dimension+1)^2+1+|W|. \qedhere
    \]
\end{proposition}
\begin{proof}
Let 
\[
  B_0 \subset (\CP^\dimension)^{\maxrank(\tcd)+1} / \operatorname{PGL}(\dimension+1)
\]
be the Zariski-open subset of point configurations satisfying the condition in Step~(2).  
Its dimension is
\begin{align*}
  \dim B_0 
  &= \dim (\CP^\dimension)^{\maxrank(\tcd)+1}
     - \dim \operatorname{PGL}(\dimension+1) \\
  &= \dimension \, (\maxrank(\tcd) + 1)
     - (\dimension + 1)^2 + 1.
\end{align*}

For $1 \le j \le |W| - \maxrank(\tcd)$, let
\[
  B_j \subset (\CP^\dimension)^{\maxrank(\tcd) + 1 + j}
\]
denote the set of point configurations obtainable after $j$ iterations of Step~(3). The construction algorithm relates $B_{j+1}$ and $B_j$ as follows.  
Let $w_1$ be the white vertex considered in the $(j+1)$-st iteration, and let 
$U_j \subset B_j$ be the open subset where 
$\tcdmap(w_2) \neq \tcdmap(w_3)$.  
Consider the $\CP^1$-bundle over $U_j$ whose fibers are the lines 
joining $\tcdmap(w_2)$ and $\tcdmap(w_3)$, and let 
$B_{j+1}$ be the open subset where 
$\tcdmap(w_1) \notin \{\tcdmap(w_2), \tcdmap(w_3)\}$.

By construction, each element of 
$B_{|W| - \maxrank(\tcd) - 1}$ corresponds to a TCD map modulo 
$\operatorname{PGL}(\dimension+1)$, and conversely every such map determines a unique element of this space. Finally, the dimension of $B_{|W| - \maxrank(\tcd) - 1}$ equals
\[
  \dim B_0 + |W| - \maxrank(\tcd) - 1,
\]
which proves the claim.
\end{proof}

In the remainder of this section, we examine how TCD maps transform under central projections. This will allow us to prove results for TCD maps of maximal rank and then obtain results for general rank by central projection.

Recall that, given a center $\subspacea \subset \CP^\dimension$, we defined the central projection $\pi_\subspacea$ in Definition~\ref{def:projection}. 
Now consider the projection $\pi_\subspacea(\tcdmap)$ of a TCD map~$\tcdmap$. 
Generically, $\pi_\subspacea(\tcdmap)$ is again a TCD map; however, this is not always the case.

\begin{definition}[$\tcdmap$-admissible subspaces]
  For a given TCD map~$\tcdmap$, we say that a projective subspace $\subspacea$ is \emph{$\tcdmap$-admissible} if
  \begin{enumerate}
    \item $\tcdmap(w) \notin \subspacea$ for all white vertices $w$.
    \item $\subspacea \cap \projectiveline_b = \varnothing$ for all black vertices $b$.\qedhere
  \end{enumerate}
\end{definition}

It is straightforward to see that $\pi_\subspacea(\tcdmap)$ is a well-defined TCD map if and only if $\subspacea$ is $\tcdmap$-admissible.  
If the first condition fails, then $\pi_\subspacea(\tcdmap(w))$ is not defined. If the second condition fails, then the images $\pi_\subspacea(\tcdmap(w))$ for the three white vertices $w$ adjacent to a black vertex~$b$ coincide, which is not allowed in Definition~\ref{def:tcdmap}.

Note that the set of $\tcdmap$-admissible subspaces forms a Zariski-open subset of the Grassmannian $\operatorname{Gr}(\dimension-r,\dimension)$.

We now show that any lower-rank TCD map arises as a central projection of a maximal-rank TCD map. This will allow us to prove statements in the maximal-rank case and obtain the general case by projecting.

\begin{proposition}\label{prop:central_projection_lift_exists}
Let $1 \leq r <  \dimension \leq \maxrank(\tcd)$. Then any TCD map $\tcdmap \in \maps(\tcd,r)$ is the central projection of a TCD map $\hat \tcdmap \in \maps(\tcd,\dimension)$.
\end{proposition}

\begin{proof}
    We may fix the center $\subspacea$ and complementary subspace $\subspaceb$ to be
    \[
       \subspacea= \{[x] \in \CP^\dimension: x_1=\cdots=x_{r+1}=0\},
      \qquad
       \subspaceb=\{[x]\in\CP^\dimension: x_{r+2}=\cdots=x_{\dimension+1}=0\}.
    \]
    Let $(\vrc,\relation)$ be a VRC lift of $\tcdmap$. We follow the steps of the construction algorithm to lift each $\vrc(w) \in \subspaceb$. In Step~(2), choose the vectors $\mathsf{u}({w_i^\partial}) \in \C^{\dimension-r}$ associated to minimal white vertices such that 
    \[
        \Span \{(\vrc(w_i^\partial), \mathsf{u}({w_i^\partial})):  w_i^\partial \in W^{\min} \} = \C^{\dimension+1}.
    \]
    In each iteration of Step~(3), suppose the relation $\relation_b$ is 
    \[
      \alpha_1 \vrc(w_1) + \alpha_2 \vrc(w_2) + \alpha_3 \vrc(w_3) = 0.
    \]
    Then, let $\mathsf{u}(w_1) := -\frac{\alpha_2}{\alpha_1} \mathsf{u}(w_2) - \frac{\alpha_3}{\alpha_1} \mathsf{u}(w_3)$.
    The conclusion follows from defining for every white vertex $w$
    \begin{equation} \label{eq:vrc_lift_central_projection}
        \hat \vrc(w):=  (\vrc(w), \mathsf{u}(w)), \qquad \hat \tcdmap(w):=[\hat \vrc(w)].\qedhere
    \end{equation}
\end{proof}

\section{Moves for TCD maps}\label{sec:tcddskp}

In Section~\ref{sec:tcds}, we defined 2-2 moves for TCDs. In this section we show how a 2-2 move $s:\tcd \rightsquigarrow \tilde{\tcd}$ at a bigon of a TCD gives rise to a move on the corresponding TCD maps defined on $\tcd$. There are two types of moves for TCD maps, depending on the orientation of the bigon of the TCD:
\begin{enumerate}
    \item The \emph{\resplitmove{}} (clockwise bigon);
    \item The \emph{spider move} (counterclockwise bigon).
\end{enumerate}
For both types of moves, we will show that the assignment $\tcdmap \mapsto \tilde \tcdmap$ defines a birational map on the corresponding moduli spaces
\[
\mu_s: \maps(\tcd,\dimension) \dashrightarrow \maps(\tilde{\tcd},\dimension),
\]
that is compatible with central projections.

\subsection{The \resplitmove{}} \label{sec:resplit}

\begin{figure}
	\centering	
	\raisebox{-0.5\height}{\includegraphics[angle=90,scale=0.8]{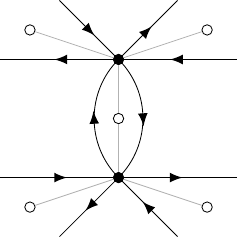}}
	\hspace{2mm}$\leftrightarrow$\hspace{1mm}
	\raisebox{-0.5\height}{\includegraphics[scale=0.8]{tikz/tcdrespliteye}}
	\hspace{15mm}
	\raisebox{-0.5\height}{\includegraphics[scale=1]{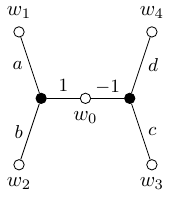}}
	\hspace{4mm}$\leftrightarrow$\hspace{3mm}
	\raisebox{-0.5\height}{\includegraphics[scale=1]{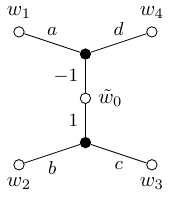}}
	\caption{The TCD $\tcd$ and graph $\pb$ as they appear in the 2-2 move at a clockwise bigon (left) and a counterclockwise bigon (right).}\label{fig:tcd_move_resplit}
\end{figure}

Consider a TCD map $\tcdmap$ on a TCD $\tcd$ with VRC $(\vrc, \relation)$. 
A \resplitmove{} corresponds to a 2-2 move on $\tcd$ at a clockwise face (see Figure~\ref{fig:tcd_move_resplit}). At the level of the VRC, it can be viewed as the composition of a contraction and a split-hence the name. The vectors of the new VRC $(\tilde\vrc,\tilde \relation)$ are
\begin{equation}
  \tilde\vrc(w_i) = \vrc(w_i) \quad (i=1,2,3,4), 
  \qquad
  \tilde\vrc(\tilde w_0)
  = a\tilde\vrc(w_1) + d\tilde\vrc(w_4)
  = -b\tilde\vrc(w_2) - c\tilde\vrc(w_3),
  \label{eq:vrc_resplit}
\end{equation}
and the new edge weights are, after a gauge transformation, as shown in Figure~\ref{fig:tcd_move_resplit}. The new TCD map $\tilde\tcdmap$ is obtained by projectivizing $\tilde\vrc$.

\begin{definition}[\resplitmovecap{} on TCD maps]\label{def:resplit-move}  
    The \emph{\resplitmove{} on TCD maps} is the rational map 
    \[
        \mu_s: \maps(\tcd,\dimension) \dashrightarrow \maps(\tilde{\tcd},\dimension),
    \]
    that sends a TCD map $\tcdmap$ to the TCD map $\tilde\tcdmap$. 
\end{definition}

The indeterminacy locus of $\mu_s$ is the closed set where $\tcdmap(w_1)=\tcdmap(w_4)$ or $\tcdmap(w_2)=\tcdmap(w_3)$, since in these cases the move produces coincident white neighbors at a trivalent black vertex and fails to define a TCD map.
The reverse resplit $s^{-1}:\tilde\tcd\rightsquigarrow\tcd$ gives a rational inverse with indeterminacy locus where $\tilde\tcdmap(w_1)=\tilde\tcdmap(w_3)$ or $\tilde\tcdmap(w_2)=\tilde\tcdmap(w_4)$.
Hence $\mu_s$ is birational.

\begin{remark}
    The indeterminacy locus of $\mu_s$ is empty whenever the rank of $\tcdmap$ is maximal. Indeed, if $\tcdmap$ has maximal rank, then the five points participating in a resplit must span a plane. The incidence relations then imply that if $\tcdmap(w_1)=\tcdmap(w_4)$, then they would also coincide with $\tcdmap(w_0)$, contradicting the definition of a TCD map. A further consequence is that any TCD map that fails to be move-generic due to a non-performable resplit cannot have maximal rank. Moreover, this shows that this type of singularity can always be resolved by lifting to a TCD map of maximal rank.
\end{remark}

\begin{lemma}\label{lem:resplit-move-projection-square}
	Let $s:\tcd\rightsquigarrow\tilde\tcd$ be a \resplitmove{} and let $\pi_{\subspacea}$ be a central projection. Then the following diagram 
	\[
		\begin{tikzcd}
		\maps(\tcd,\dimension\!+\!1) \arrow[r, dashed, "\mu_s"] \arrow[d, dashed, "\pi_{\subspacea}"']
		& \maps(\tilde\tcd,\dimension\!+\!1) \arrow[d, dashed, "\pi_{\subspacea}"] \\
		\maps(\tcd,\dimension) \arrow[r, dashed, "\mu_s"']
		& \maps(\tilde\tcd,\dimension)
		\end{tikzcd}
	\]
	commutes.
\end{lemma}

\begin{proof}
    The definition of $\tilde \vrc$ in \eqref{eq:vrc_resplit} is compatible with the quotient map \eqref{eq:central_projection_affine}.
\end{proof}

It turns out the resplit is closely related to an equation which plays a prominent role in discrete integrable systems.

\begin{definition}[Discrete Schwarzian Kadomtsev-Petviashvili equation]\label{def:dskpequation}
    Let $\point_1, \point_{1,2}, \point_2, \point_{2,3},\point_3,\point_{3,1}$ be $6$ points in $\CP^1$. We say that the points satisfy the \emph{discrete Schwarzian Kadomtsev-Petviashvili equation} or \emph{dSKP equation} \cite{ncwqdskp, dndskp, ksclifford} for short if 
    \[
        \mr(\point_1, \point_{1,2}, \point_2, \point_{2,3},\point_3,\point_{3,1})=-1. \qedhere
    \]
\end{definition}
The relation to TCD maps is due to the following proposition.
\begin{proposition}\label{prop:resplitmr}
	Let $\tcdmap(w_0),\tilde \tcdmap(\tilde w_0),\tcdmap(w_1),\tcdmap(w_2),\tcdmap(w_3),\tcdmap(w_4)$ be the $6$ points involved in a resplit as in Figure~\ref{fig:tcd_move_resplit} (left). Then 
	\[
		\mr(\tcdmap(w_1),\tcdmap(w_0),\tcdmap(w_2),\tcdmap(w_3),\tilde \tcdmap(\tilde w_0),\tcdmap(w_4)) = -1. \label{eq:resplitmr}  \qedhere
	\]
\end{proposition}
\proof{
    This follows by a direct calculation using~\eqref{eq:vrc_resplit}.
}

If the rank of $\tcdmap$ is 1, then Proposition~\ref{prop:resplitmr} states that the 6 points satisfy the dSKP equation. Moreover, if $\tcdmap$ has higher rank, one can always obtain central projections of rank 1, which can be interpreted as coordinate projections of $\tcdmap$. In particular, by Lemma~\ref{lem:projinvariants} (3), the coordinates of a TCD map under the resplit satisfy the dSKP equation. 

\begin{remark} 
    In the rank 1 case, a geometric interpretation of the dSKP equation based on Clifford’s four-circle theorem was given in \cite{ksclifford}, under which the resplit move can be described in terms of intersections of certain circles. 
\end{remark}

\begin{figure}
	\includegraphics{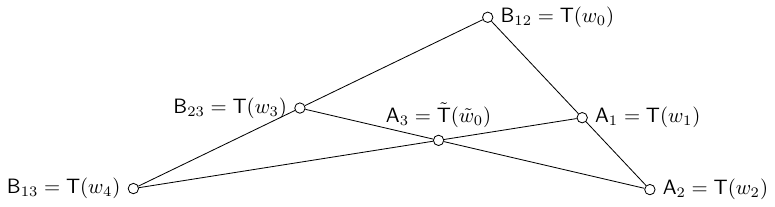}
	\caption{A Menelaus configuration associated to a \resplitmove{} in the case that the points span a two-dimensional space.}
	\label{fig:menelaus}
\end{figure}

If the TCD map is not rank 1, then the points involved in a resplit may span a plane. In this case, the new point $\tilde \tcdmap(\tilde w_0)$ can be defined as an intersection of lines. Specifically,
\begin{align*}
    \tilde \tcdmap(\tilde w_0) = (\tcdmap(w_1)\tcdmap(w_4) ) \cap (\tcdmap(w_2)\tcdmap(w_3) ).
\end{align*}
In this case, Proposition~\ref{prop:resplitmr} is a consequence of the classical \emph{Menelaus theorem}; see Figure~\ref{fig:menelaus} and~\cite{ksclifford}.

\subsection{The spider move.}

\begin{figure}
	\centering
	\raisebox{-0.5\height}{\includegraphics[angle=90,scale=0.8]{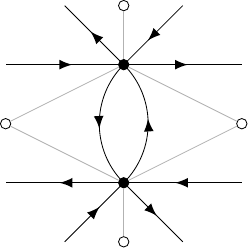}}
	\hspace{2mm}$\leftrightarrow$\hspace{1mm}
	\raisebox{-0.5\height}{\includegraphics[scale=0.8]{tikz/tcdspidereye}}
	\hspace{7mm}
	\raisebox{-0.5\height}{\includegraphics[scale=1]{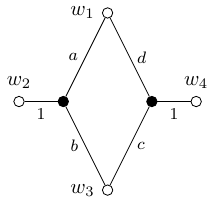}}
	\hspace{0mm}$\leftrightarrow$\hspace{0mm}
	\raisebox{-0.5\height}{\includegraphics[scale=1]{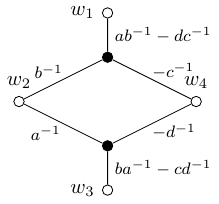}}
	\vspace{-2mm}
	\caption{Transformation of edge weights in the resplit (left) and spider move (right).}\label{fig:tcd_move_spider}
\end{figure}

The spider move corresponds to a 2-2 move on $\tcd$ at a counterclockwise bigon (see Figure~\ref{fig:tcd_move_spider}). 
At the level of the VRC, this is precisely the spider move defined earlier. 
Projectivizing the resulting VRC gives the corresponding transformation of TCD maps. 
The vectors remain unchanged, and the new edge weights are as in Figure~\ref{fig:tcd_move_spider}, after a gauge transformation. 
Note that all points of the TCD map involved in the spider move lie on a projective line.

\begin{definition}[Spider move on TCD maps]\label{def:spider-move}
    The \emph{spider move on TCD maps} is the rational map
    \[
       \mu_s:\ \maps(\tcd,\dimension)\dashrightarrow \maps(\tilde\tcd,\dimension),
    \]
    that sends $\tcdmap$ to $\tilde \tcdmap$.
\end{definition}
The indeterminacy locus of $\mu_s$ is the closed set where $\tcdmap(w_2)=\tcdmap(w_4)$ (the move would create coincident white neighbors at a trivalent black vertex). 
The inverse spider move $s^{-1}:\tilde\tcd\rightsquigarrow\tcd$ defines a rational inverse with indeterminacy locus where $\tilde\tcdmap(w_1)=\tilde\tcdmap(w_3)$. 
Hence, $\mu_s$ is birational.

\begin{remark}
    Unlike in the case of the resplit, the indeterminacy locus of $\mu_s$ is not necessarily empty for the spider move even if $\tcdmap$ has maximal rank. In fact, it follows from our definition of the projective $X$-variables (see Section~\ref{sec:projcluster}) that the indeterminacy locus of the spider move coincides with the indeterminacy locus of cluster mutation.    
\end{remark}

\begin{remark}
    As mentioned, for the spider move the points of $\tcdmap$ are unchanged and neither do new lines appear; only the combinatorics of the TCD changes. For this reason, we interpret the spider move as a \emph{reparametrization}: the underlying geometric object remains the same, but its parametrization changes. 
\end{remark}

Since the spider move does not change the points, we trivially get:
\begin{lemma}\label{lem:spider-move-projection-square}
    Let $s:\tcd\rightsquigarrow\tilde\tcd$ be a spider move and let $\pi_{\subspacea}$ be a central projection.
    Then the following diagram  
    \[
        \begin{tikzcd}
            \maps(\tcd,\dimension\!+\!1) \arrow[r, dashed, "\mu_s"] \arrow[d, dashed, "\pi_\subspacea"']
            & \maps(\tilde\tcd,\dimension\!+\!1) \arrow[d, dashed, "\pi_\subspacea"] \\
            \maps(\tcd,\dimension) \arrow[r, dashed, "\mu_s"']
            & \maps(\tilde\tcd,\dimension)
        \end{tikzcd}
    \]
    commutes.
\end{lemma}

\subsection{Move graph}

\begin{definition}[Move graph]
    Let $\complex_1$ be the graph whose vertices are minimal TCDs with a fixed strand permutation, and whose edges correspond to $2$-$2$ moves. We refer to $\complex_1$ as the \emph{move graph}.
\end{definition}

Since performing a $2$-$2$ move twice returns to the original TCD, we include only one undirected edge for both directions of the move. 

\begin{remark} \label{rem:OPS}
    The move graph of a minimal TCD is finite. Its vertices were described explicitly by Oh--Postnikov--Speyer~\cite{OPS} (see Theorem~\ref{thm:independentsets}).
\end{remark}

In some contexts, we wish to assume that moves can be performed freely, meaning that we may assume that no singularities occur.

\begin{definition}\label{def:gentcdmap}
	Let $\tcd$ be a minimal TCD and $\tcdmap: \tcd \rightarrow \CP^\dimension$ be a TCD map. We call $\tcdmap$ \emph{move-generic} if every possible sequence of 2-2 moves starting at $\tcdmap$ is well-defined.
\end{definition}

\begin{remark}
    It is not obvious that the space of move-generic TCD maps is nonempty since there are infinitely many paths in the move graph, and so the definition involves an infinite intersection of Zariski-open subsets. The non-emptiness of this space will be established in Section~\ref{sec:tcdconsistency} as a consequence of multi-dimensional consistency.
\end{remark}

\section{Multi-dimensional consistency, Desargues maps and dSKP lattices}
\label{sec:tcdconsistency}

In this section we place TCD maps within the combinatorial geometry of the \(A_{n-1}\) lattice and use this viewpoint to prove multi-dimensional consistency. Finally, we relate TCD maps to two classical lattice geometries: Desargues maps and dSKP lattices.

\subsection{The $A_n$ lattice}\label{sec:alattice} \label{sec:abel}

Let
\begin{align*}
	A_{n-1} :=  \left\{(z_1,z_2,\dots,z_{n})\in\Z^{n}\ : \ \sum_{i=1}^{n}z_i = 0\right\}
\end{align*}
denote the root lattice of type $A_{n-1}$. Replacing the coordinate sum equals $0$ condition with another integer $j$, we obtain a shifted copy of $A_{n-1}$ which we denote $A_{n-1}^j$. 

Let us explain a natural way to embed the black vertices, white vertices and faces of $\pb$ associated to any TCD into shifted copies of the $A_{n-1}$ lattice, where $n$ is the number of strands in the TCD, defined as follows. Let $\pb$ be a minimal \plabic{} graph. Recall that the faces of the medial graph $\pb^\times$ are in bijection with $B \sqcup W \sqcup F$. Define $\an: B\sqcup W\sqcup F \rightarrow \Z^n$ by
\begin{align}
	\an(x)_i =
	\begin{cases}
		1 & \text{if $x$ lies to the left of medial strand $i$},\\
		0 & \text{if $x$ lies to the right of medial strand $i$},
	\end{cases}
	\label{eq:shiftamap}
\end{align}
for every vertex and face of $\pb$. Let $k:=\maxrank(\tcd)+1$. An immediate consequence is \cite[Proposition~4.2]{mstwist},
\[
    \an(B) \subset A_{n-1}^{k+1},  \qquad  \an(W) \subset A_{n-1}^{k-1},  \qquad  \an(F) \subset A_{n-1}^{k}.
\]

\begin{figure}
	\centering
	\raisebox{-0.5\height}{\includegraphics[scale=0.8]{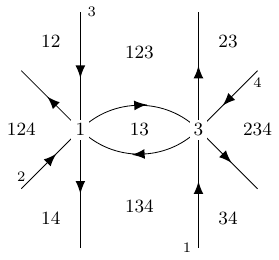}}
	\hspace{0mm}$\leftrightarrow$\hspace{0mm}
	\raisebox{-0.5\height}{\includegraphics[scale=0.8]{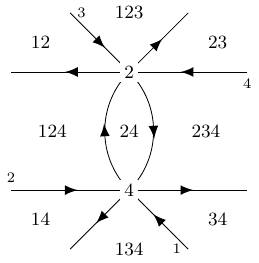}}
	\hspace{6mm}
	\raisebox{-0.5\height}{\includegraphics[scale=0.8]{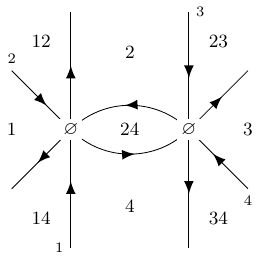}}
	\hspace{0mm}$\leftrightarrow$\hspace{0mm}
	\raisebox{-0.5\height}{\includegraphics[scale=0.8]{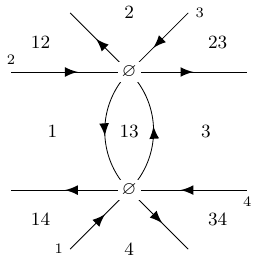}}

	\caption{The transformation of labels of clockwise faces of $\tcd$ (white vertices of $\pb$) in the \resplitmove{} (left) and the spider move (right), together with the labeling $\ban$. 
    } 	\label{fig:an_transform}
\end{figure}

However, for our purposes it will be more practical to work with a slightly modified version of $\an$ which is defined as follows:
\begin{align*}
    \ban(x)_i = 
        \begin{cases}
            \an(x)_i & \mbox{if } x \in W \sqcup F, \\
            \an(x)_i  & \mbox{if } x \in B \mbox{ and the strand $i$ is not incident to $x$},   \\
            0  & \mbox{if } x \in B \mbox{ and the strand $i$ is incident to $x$}.
        \end{cases}
\end{align*}

As a consequence,
\[
    \ban(B) \subset A_{n-1}^{k-2},  \qquad  \ban(W) \subset A_{n-1}^{k-1},  \qquad  \ban(F) \subset A_{n-1}^{k}.
\]

\begin{remark}
    Equivalently, we can assign to each vertex or face of $\pb$ a subset of $\{1,\dots,n\}$. This labeling was introduced in the study of cluster algebras and Grassmannians by Scott~\cite{scottgrass}. A version of $\an$ for \plabic{} graphs/TCDs on a torus was introduced by Fock~\cite{fockinverse} under the name \emph{discrete Abel map}. Versions of the map $\ban$ for planar graphs appear in the inverse boundary measurement map \cite[(5)]{GonKon} and for toric graphs in the inverse spectral transform (see~\cite[Section~4]{fockinverse}~and~\cite[Section 5.1]{GGK}).
\end{remark}

The map $\ban$ is compatible with the 2-2 moves in the sense that $\ban$ only changes locally (see Figure~\ref{fig:an_transform}). Specifically, for the resplit the restriction $\ban|_F$ is unchanged; $\ban|_W$ only changes for the white vertex of the vector that is replaced in the resplit; and $\ban|_B$ only changes for the two lines that are replaced. For the spider move there is no change to $\ban|_W$ and $\ban|_B$ reflecting the fact that the geometry is unchanged, and $\ban|_F$ only changes for the center face. We will see in Section~\ref{sec:clusterstructures} that this reflects that the $X$-variable of the center face changes when performing the spider move.

\subsection{Cycles in the move graph}

Starting from a given TCD, one can perform a sequence of 2--2 moves (without reversing any move) and return to the original TCD. In other words, the move graph $\complex_1$ contains non-trivial cycles. It was conjectured by Thurston~\cite{thurstontriple} that these cycles are generated by the following.

\begin{figure}
	\centering
	\includegraphics[scale=0.85]{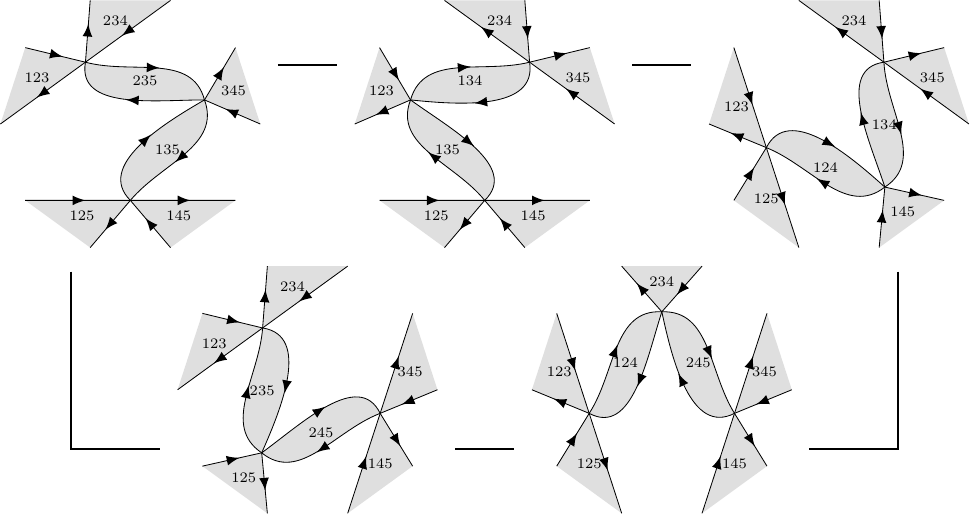}
	
	\caption{An elementary 5-cycle in the move graph of TCDs with strand permutation $\enm{5}{4}$. Reversing the orientations of strands gives the other elementary 5-cycle with strand permutation $\enm{5}{2}$.
    }
	\label{fig:tcdmovegraphfivecycle}
\end{figure}

\begin{figure}
	\centering
	
	\includegraphics{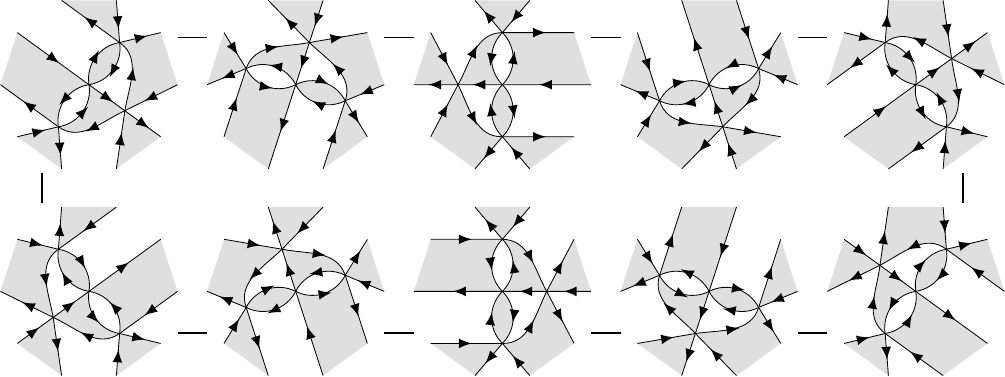}	
	\caption{An elementary 10-cycle in the move graph of TCDs with strand permutation $\enm{5}{3}$.} 
	\label{fig:tcdmovegraphtencycle}
\end{figure}

\begin{definition}[Elementary cycles]
We call the following cycles \emph{elementary}:
\begin{enumerate}
    \item Every $4$-cycle arising in a subdiagram that is a TCD containing two disjoint copies of one of the configurations in Figure~\ref{fig:an_transform}, where two 2--2 moves can be performed in either order.
    \item Every $5$-cycle arising in a subdiagram that is a TCD with strand permutation $\enm{5}{4}$ (Figure~\ref{fig:tcdmovegraphfivecycle}) or $\enm{5}{2}$ (the same diagram with reversed strand orientations).
    \item Every $10$-cycle arising in a subdiagram that is a TCD with strand permutation $\enm{5}{3}$ (see Figure~\ref{fig:tcdmovegraphtencycle}).\qedhere
\end{enumerate}
\end{definition}

\begin{definition}[Move graph 2-complex]
Let $\complex_2$ denote the $2$-complex obtained from $\complex_1$ by attaching:
\begin{itemize}
    \item A square along each elementary $4$-cycle.
    \item A pentagon along each elementary $5$-cycle.
    \item A decagon along each elementary $10$-cycle.\qedhere
\end{itemize}
\end{definition}

Thurston's conjecture was proved by Balitskiy and Wellman~\cite{bwtriple}:

\begin{theorem}[\cite{bwtriple}]\label{th:tcdtwocycles}
    The complex $\complex_2$ is simply connected.
\end{theorem}

 Multi-dimensional consistency is a flatness condition on TCD maps over the move-graph and is a hallmark of discrete integrability; see also~\cite{absintquads,absquads, ddgbook}.

 \begin{definition}[Multi-dimensional consistency]
    We say that a TCD map is \emph{(multi-dimensionally) consistent} if the composition of local moves for TCDs around any closed cycle in the move graph yields the initial TCD map.
\end{definition}

There are many geometric discrete dynamical systems that are defined by starting with two-dimensional boundary data and applying local propagation rules to obtain $n$-dimensional solutions; for example, the propagation of Q-nets~\cite{doliwasantiniqnet}, Darboux maps~\cite{schieflattice} and line complexes~\cite{bobenkoschieflinecomplexes}. However, when $n>3$, there are typically many ways to propagate from initial data to the full data. Hence, one must verify that the result of propagation does not depend on the chosen sequence. The system is multi-dimensionally consistent if the result of propagation is independent of the choice of sequence for all choices of initial data.

As we explain in~\cite{paper2}, all known multi-dimensionally consistent geometric discrete dynamical systems with two-dimensional initial data can be realized as special cases of TCD maps. Moreover, as shown in Theorem~\ref{th:tcd_consistency_along_fund_cycles}, TCD maps themselves are multi-dimensionally consistent. This provides a unified geometric proof of multi-dimensional consistency for the aforementioned systems.

Theorem~\ref{th:tcdtwocycles} allows us to give an elementary proof of the multi-dimensional consistency of TCD maps. Additionally, the following classical theorem will play a key role.
 
 \begin{theorem}[Desargues' theorem]\label{th:desargues}
Let $\mathsf{A}_1, \mathsf{A}_2, \mathsf{A}_3$ and $\mathsf{B}_1, \mathsf{B}_2, \mathsf{B}_3$ be the vertices of two triangles in $\CP^\dimension$ with $\dimension \geq 2$.  
The following statements are equivalent:
\begin{enumerate}
	\item The lines $\mathsf{A}_1\mathsf{B}_1$, $\mathsf{A}_2\mathsf{B}_2$ and $\mathsf{A}_3\mathsf{B}_3$ meet at a point.
	\item The points $\mathsf{A}_1\mathsf{A}_2 \cap \mathsf{B}_1\mathsf{B}_2$,  
	      $\mathsf{A}_2\mathsf{A}_3 \cap \mathsf{B}_2\mathsf{B}_3$ and  
	      $\mathsf{A}_3\mathsf{A}_1 \cap \mathsf{B}_3\mathsf{B}_1$ lie on a line.\qedhere
\end{enumerate}
\end{theorem}

\begin{figure}
	\includegraphics{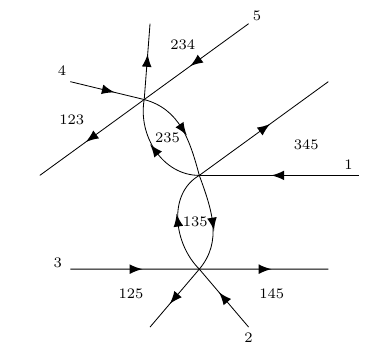}
	\hspace{.2cm}
	\includegraphics{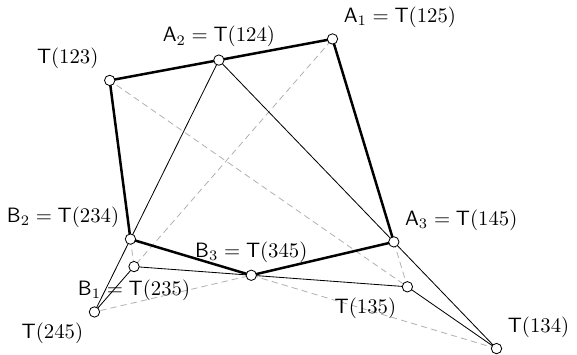}

	\caption{Labeling of strands and clockwise faces in the top TCD in Figure~\ref{fig:tcdmovegraphfivecycle}~(left) and a Desargues configuration (right).}
	\label{fig:desargues}
\end{figure}

\begin{lemma} \label{lem:tcd_consistency_along_fund_cycles}
	Move-generic TCD maps are consistent along each elementary cycle.
\end{lemma}
\begin{proof}

    Consistency around $4$-cycles is clear since the moves act on disjoint parts of the diagram. 
    
    In the cases of the $5$- and $10$-cycles, it suffices to verify consistency on the corresponding subdiagrams, whose move graph is exactly one of the elementary cycles. Moreover, it suffices to prove the maximal-rank case; for lower-rank cases, first take a maximal-rank lift provided by Proposition~\ref{prop:central_projection_lift_exists}, apply the maximal-rank case and then descend using Lemma~\ref{lem:resplit-move-projection-square} together with the invariance of multi-ratios under central projections (Lemma~\ref{lem:projinvariants}).

    There are two 5-cycle configurations. The $5$-cycle of TCDs with strand permutation $\enm 52$ has no internal white vertices and is therefore trivially consistent. Thus the interesting case is the $\enm 54$ strand permutation. By Proposition~\ref{prop:max_rank}, the maximal rank of any TCD map in this cycle is $7-3-1=3$ and it gives rise to a Desargues configuration as follows. Let $\tcd$ be the TCD shown in Figure~\ref{fig:desargues} (left). This is the top left TCD in Figure~\ref{fig:tcdmovegraphfivecycle}. Consider the corresponding configuration of points and lines in Figure~\ref{fig:desargues} (right), 
    with each object labeled by $\ban$. Now apply 2-2 moves at the clockwise faces (\resplitmove{}s) as follows to move three steps counterclockwise around Figure~\ref{fig:tcdmovegraphfivecycle} to reach the top right TCD:  
    \begin{itemize}
        \item $135$: create the new point $\tcdmap(245) := \tcdmap(125) \tcdmap(235) \cap \tcdmap(145)\tcdmap(345)$.
        \item $235$: create the new point $\tcdmap(124) := \tcdmap(123)\tcdmap(125) \cap \tcdmap(234) \tcdmap(245)$.
        \item $245$: create the new point $\tcdmap(134) := \tcdmap(124)\tcdmap(145) \cap \tcdmap(234) \tcdmap(345)$. 
    \end{itemize}
    
    Now let
    \[
    \begin{alignedat}{3}
    \mathsf{A}_1 &:= \tcdmap(125), &\qquad \mathsf{A}_2 &:= \tcdmap(124), &\qquad \mathsf{A}_3 &:= \tcdmap(145), \\
    \mathsf{B}_1 &:= \tcdmap(235), &\qquad \mathsf{B}_2 &:= \tcdmap(234), &\qquad \mathsf{B}_3 &:= \tcdmap(345).
    \end{alignedat}
    \]
    Define the lines between the corresponding vertices by
    \[
    \projectiveline(25):=\mathsf{A}_1\mathsf{B}_1,\qquad
    \projectiveline(24):=\mathsf{A}_2\mathsf{B}_2,\qquad
    \projectiveline(45):=\mathsf{A}_3\mathsf{B}_3,
    \]
    which meet at the common point $\tcdmap(245)$. Denote the three sides of the $\mathsf A$-triangle and of the $\mathsf B$-triangle by
    \[
    \projectiveline(12):=\mathsf{A}_1\mathsf{A}_2,\qquad
    \projectiveline(14):=\mathsf{A}_2\mathsf{A}_3,\qquad
    \projectiveline(15):=\mathsf{A}_3\mathsf{A}_1,
    \]
    \[
    \projectiveline(23):=\mathsf{B}_1\mathsf{B}_2,\qquad
    \projectiveline(34):=\mathsf{B}_2\mathsf{B}_3,\qquad
    \projectiveline(35):=\mathsf{B}_3\mathsf{B}_1.
    \]
    With this notation, the intersections of corresponding sides are
    \[
    \mathsf A_1\mathsf A_2 \cap \mathsf B_1\mathsf B_2
    =\projectiveline(12)\cap\projectiveline(23)=\tcdmap(123),
    \]
    \[
    \mathsf A_2\mathsf A_3 \cap \mathsf B_2\mathsf B_3
    =\projectiveline(14)\cap\projectiveline(34)=\tcdmap(134),
    \]
    \[
    \mathsf A_3\mathsf A_1 \cap\mathsf B_3\mathsf B_1
    =\projectiveline(15)\cap\projectiveline(35)=\tcdmap(135).
    \]
    By Desargues’ theorem (Theorem~\ref{th:desargues}), these three points lie on a line, which we denote by $\projectiveline(13)$. Figure~\ref{fig:desargues}~(right) shows the resulting configuration of points and lines. 
    
    On the other hand, we can apply 2-2 moves as follows moving two steps clockwise around Figure~\ref{fig:tcdmovegraphfivecycle}, again arriving at the top right TCD:
    \begin{itemize}
    \item 235: create the new point $\mathsf{X}:= \mathsf{T}(123)\mathsf{T}(135) \cap \mathsf{T}(234) \mathsf{T}(345)$ which is equal to $\projectiveline(13) \cap \projectiveline(34) = \mathsf{T}(134)$.
    \item 135: create the new point $\mathsf{Y}:= \mathsf{T}(123)\mathsf{T}(125) \cap \mathsf{T}(145) \mathsf{X}$ which is equal to $\projectiveline(12) \cap \projectiveline(14) = \mathsf{T}(124)$.
    \end{itemize}
    Thus both paths give the same TCD map on the top right TCD, proving consistency.

     By Proposition~\ref{prop:max_rank}, the maximal rank of any TCD map in the 10-cycle is $7-4-1=2$. The points attached to the two internal white vertices of any TCD map in the cycle are uniquely determined by the points attached to the boundary vertices by incidence conditions. Since the boundary points do not change under 2-2 moves, it is consistent around the cycle.
\end{proof}

\begin{theorem} \label{th:tcd_consistency_along_fund_cycles}
    Move-generic TCD maps on a TCD $\tcd$ are consistent. Moreover, they form a non-empty Zariski-open subset of the space of TCD maps on $\tcd$. 
\end{theorem}

\begin{proof}
    Recall from Remark~\ref{rem:OPS} that the move graph $\complex_1$ of $\tcd$ is finite. Fix a spanning tree $\mathfrak T$ of $\complex_1$ rooted at $\tcd$. For each vertex $\tcd'$, there is a unique path along $\mathfrak T$ from $\tcd$ to $\tcd'$. The set of TCD maps on $\tcd$ such that all the TCD maps defined along these paths are valid forms a nonempty Zariski-open subset of the space of TCD maps on $\tcd$.
    
    Considering an arbitrary path in $\complex_1$ from $\tcd$ to $\tcd'$ does not lead to the appearance of any additional TCD maps. Indeed, the cycle in the move graph obtained by concatenating the latter path with the unique path in $\mathfrak T$ from $\tcd'$ to $\tcd$ forms a cycle in $\complex_1$, which can be decomposed into finitely many elementary cycles by Theorem~\ref{th:tcdtwocycles}. Lemma~\ref{lem:tcd_consistency_along_fund_cycles} then guarantees consistency along each elementary cycle, thus along the whole cycle.
\end{proof}

\begin{remark}
  By Theorem~\ref{th:tcd_consistency_along_fund_cycles}, it is natural and well-defined to label the points and lines of a TCD map by $\ban$ like we did in the proof of Lemma~\ref{lem:tcd_consistency_along_fund_cycles}, and this labeling remains compatible under moves.
\end{remark}

\begin{remark}
    It is not surprising that consistency is governed by Desargues' theorem. For instance, in the case of fundamental line complexes, it was already noted that multi-dimensional consistency is due to Desargues' theorem in \cite{bobenkoschieflinecomplexes}. In the next section, we will study Desargues maps, which are so named because their consistency also follows from Desargues’ theorem. Moreover, by considering TCD maps in $\CP^1$, we obtain a new combinatorial proof of the multi-dimensional consistency of the dSKP equation. A (different) proof of this was first given by Adler, Bobenko, and Suris in their classification of integrable equations of octahedral type~\cite{absoctahedron}. 
Later, King and Schief showed that the consistency of the dSKP equation follows geometrically from a conformal version of Desargues’ theorem~\cite{ksconformaldesargues}, which fits naturally with the interpretation that we give here.
\end{remark}

\subsection{Desargues maps}

Consider a TCD $\tcd$ with corresponding BTB graph $\pb$ and let $k= \maxrank(\tcd)+1$. Recall that, via the map $\ban$, we identify the set of white vertices of $\pb$ with a subset of $A_{n-1}^{k-1}$.
The 2-cells of the $A_{n-1}^{k-1}$ lattice can be partitioned into two families usually referred to as white triangles and black triangles. In this section, we need the \emph{black triangles}, which are given by three lattice points of the form
\[
    \{ z + e_i : i \in I \} \subset A_{n-1}^{k-1},
\]
where $z \in A_{n-1}^{k-2} $ and $ I \subset \{1,2,\dots,n\}$ is a $3$-element subset.

\begin{definition}[Desargues map] \label{def:desarguesmaps}
   Let $L\subset A_{n-1}$. A \emph{Desargues map} is a map \[\mathsf{D}: L \longrightarrow \CP^\dimension,\] such that for every black triangle in $L$, the images of its three vertices under $\mathsf{D}$ lie on a common line. 
\end{definition}

Desargues maps were originally introduced by Doliwa on $\Z^{n-1}$ in a different form~\cite{doliwadesargues}. The equivalence with Definition~\ref{def:desarguesmaps} was established in~\cite[Proposition~3.1]{doliwadesarguesweyl}.
The goal of this section is to establish a relationship between TCD maps and Desargues maps. In order to do this, we need to introduce some combinatorial notions that will only be used in the proofs of Lemmas~\ref{lem:tcdthreestar} and \ref{lem:whitestar}.
The following notion was introduced by Leclerc and Zelevinsky~\cite{LZqcomm}. 
\begin{definition}[Weak separation]
    We say that $z,z' \in \Z^n$ are \emph{weakly separated} if there does not exist four indices $i_1 < i_2 < i_3 < i_4$, such that
\[
    z_{i_1} < z'_{i_1}, \quad z_{i_2} > z'_{i_2}, \quad z_{i_3} < z'_{i_3}, \quad z_{i_4} > z'_{i_4}. \qedhere
\]
\end{definition}

\begin{theorem}[Oh, Postnikov and Speyer~\cite{OPS}]\label{thm:independentsets}
Let $\tcd$ be a minimal TCD with BTB graph $\pb$, and let $L \subset A_{n-1}^{k-1}$ be the set of $\ban$-labels of all white vertices that appear in all BTB graphs move-equivalent to $\pb$. Then every pairwise weakly separated collection in $L$ that is maximal by inclusion is also maximal by size. Moreover, such a collection is the $\ban$-image of the white vertices of a BTB graph $\tilde \pb$, unique up to spider moves, move-equivalent to $\pb$, and conversely.
\end{theorem}
\begin{proof}
In~\cite{OPS}, the result is formulated in terms of face labels; however, it translates directly to white vertex labels, since, after reversing the orientation of all strands, faces and white vertices play equivalent roles in the context of TCDs. 
\end{proof}

We also require the explicit construction of the graph $\tilde \pb$ appearing in Theorem~\ref{thm:independentsets} due to~\cite{OPS} which we now describe. However, we rephrase their construction using white vertex labels. Let $\mathcal D$ be a maximal by inclusion pairwise weakly separated collection in $L$. For $z \in \mathcal D$, let $S(z)$ be the corresponding $(k-1)$-element subset of $\{1,2,\dots,n\}$. For $K$ a $k$-element subset of $\{1,2,\dots,n\}$, define the corresponding \emph{face clique}
\[
\mathcal F(K):=\{ z \in \mathcal D: S(z) \subset K\}.
\]
Similarly, for $K$ a $(k-2)$-element subset, define the corresponding \emph{black clique}
\[
\mathcal B(K):=\{z \in \mathcal D: S(z) \supset K\}.
\]
A clique is called \emph{non-trivial} if it contains at least $3$ elements. 
Note that $\mathcal F(K)$ is of the form 
\[
\{ z_K-e_{i_1},\dots,z_K-e_{i_r} \},
\]
where $z_K:= \sum_{i \in K} e_i$ and $i_1 < \cdots < i_r$. Define its \emph{boundary} 
\[
\partial \mathcal F(K) := \{ (z_1, z_2), \dots, (z_r, z_1)\}, \qquad z_j:=z_K-e_{i_j}.
\]
Similarly, every black clique $\mathcal B(K)$ is of the form 
\[
\{ z_K+e_{i_1},\dots,z_K+e_{i_r} \}
\]
for some $i_1 < \cdots < i_r$. Define its \emph{boundary} 
\[
\partial \mathcal B(K) := \{ (z_1, z_2), \dots, (z_r, z_1)\}, \qquad z_j:=z_K+e_{i_j}.
\]

Let $\bm v_1,\dots,\bm v_n \in \R^2$ be the vertices of a convex $n$-gon and define the point $\bm v(z):= \sum_{i=1}^n z_i \bm v_i$. The \emph{plabic tiling} corresponding to $\mathcal D$ is the $2$-dimensional complex $\Sigma(\mathcal D) \subset \R^2$ such that:
\begin{itemize}
    \item The vertices of $\Sigma(\mathcal D)$ are $\{ \bm v(z): z \in \mathcal D\}$;
    \item The edges of $\Sigma(\mathcal D)$ are of the form $\bm v(z) \bm v(z')$ where $(z,z')$ belongs to the boundary of a non-trivial clique;
    \item The faces of $\Sigma(\mathcal D)$ are polygons whose vertices form a non-trivial clique.
\end{itemize}
The BTB graph $\tilde \pb$ is obtained from $\Sigma(\mathcal D)$ as follows. Triangulate the faces of $\Sigma(\mathcal D)$ that correspond to black cliques. Place a black vertex inside each triangular black face and a white vertex at each vertex of $\Sigma(\mathcal D)$. Connect each black vertex to every white vertex lying on the boundary of the corresponding triangular face. Different choices of triangulations give rise to different $\tilde \pb$ that are related by spider moves.

\begin{lemma}\label{lem:tcdthreestar}
Let $\tcd$ be a minimal TCD with BTB graph $\pb$. Suppose that for some $z \in A^{k-2}_{n-1}$ and $\{i_1,i_2,i_3\} \subset \{1,2,\dots,n\}$,
\[
    z + e_{i_1}, \quad z + e_{i_2}, \quad z + e_{i_3} \in L.
\]
Then there exists a BTB graph $\tilde{\pb}$ move-equivalent to $\pb$ containing a black vertex $b$ incident to white vertices $w_1,w_2,w_3$ such that
\[
    \ban(b) = z, \qquad 
    \ban(w_j) = z + e_{i_j}, \qquad j = 1,2,3. \qedhere
\]
\end{lemma}

\begin{proof}
Since the collection $\{z+e_{i_1},z+e_{i_2},z+e_{i_3}\}$ is weakly separated, it can be extended to a maximal weakly separated collection which corresponds to a BTB graph $\tilde \pb$ by Theorem~\ref{thm:independentsets}. Moreover, $\{z+e_{i_1},z+e_{i_2},z+e_{i_3}\} \subset \mathcal B(K)$ where $K:=\{1 \leq i \leq n: z_i =1\}$. Choose a triangulation of black cliques such that they form the vertices of a triangle and let $b$ be the corresponding black vertex.
\end{proof}

With these combinatorial preliminaries in place, we can now state and prove a close connection between TCD maps and Desargues maps.

\begin{theorem}\label{th:tcddesargues}
Let $\tcd$ be a minimal TCD with BTB graph $\pb$, and let 
$L \subset A_{n-1}^{k-1}$ denote the set of $\ban$-labels of all white vertices appearing in BTB graphs that are move-equivalent to $\pb$. Then every TCD map 
\[
\tcdmap : W \longrightarrow \CP^\dimension
\]
on a \plabic{} graph $\pb$ with $n$ strands is the restriction of a Desargues map
\[
\mathsf{D} : L \longrightarrow \CP^\dimension
\]
to the subset $\ban(W) \subset L$. 
Moreover, if $\tilde{\tcdmap}$ is a TCD map on a BTB graph $\tilde{\pb}$ 
obtained from $\tcdmap$ by a sequence of local moves, then 
$\tilde{\tcdmap}$ is the restriction of the \emph{same} Desargues map $\mathsf{D}$ 
to $\ban(\tilde{W})$, where $\tilde{W}$ is the set of white vertices of $\tilde{\pb}$.
\end{theorem}
\begin{proof}
    By Lemma~\ref{lem:tcdthreestar}, if three vertices of a black triangle appear in $L$, then there is a move-equivalent TCD $\tilde \tcd$ such that they appear as white vertices sharing a black vertex. Hence, the corresponding points of $\tilde \tcdmap$ are on a line, which shows that $\tcdmap$ is the restriction of a Desargues map.
 
    The second statement follows because Desargues maps (and TCD maps) are multi-dimensionally consistent \cite{doliwadesargues}.
\end{proof}

A consequence of Theorem~\ref{th:tcddesargues} is that TCD maps can be interpreted as initial data for a Desargues map. The global propagation of a Desargues map from this initial data is governed by sequences of resplits and spider moves.

\subsection{dSKP lattices}

Recall from Section~\ref{sec:resplit} that the resplit for rank 1 TCD maps is governed by the dSKP equation. Also recall that, for a TCD map of rank $\dimension>1$, any central projection $\pi_\subspacea(\tcdmap)$ from a center $\subspacea$ of dimension $(\dimension-2)$ is a rank 1 TCD map, which can be viewed as a coordinate projection of $\tcdmap$. In this section, we explain how rank 1 TCD maps are related to \emph{dSKP lattices}~\cite{ksclifford}.

For this construction, we need the 3-cells of the $A_{n-1}^{k-1}$ lattice that are \emph{octahedra}. An octahedron is given by the six lattice points of the form
\[
    \{z + e_i + e_j : i, j \in I \mbox{ and } i < j \} \subset A_{n-1}^{k-1},
\]
where $z \in A_{n-1}^{k-3} $ and $ I \subset \{1,2,\dots,n\}$ is a $4$-element subset.

\begin{definition}[dSKP lattice]
    Let $L\subset A_{n-1}^{k-1}$. A \emph{dSKP lattice} is a map \[\mathsf S: L \longrightarrow \CP^1\] such that, for every octahedron in $L$, the images of its six vertices satisfy the dSKP equation~\eqref{eq:resplitmr}. 
\end{definition}

Instead of Lemma~\ref{lem:tcdthreestar}, we have:

\begin{lemma}\label{lem:whitestar}
Let $\tcd$ be a minimal TCD with BTB graph $\pb$. 
Suppose that for some $z \in A^{k-2}_{n-1}$ and distinct indices 
$\{i_1,i_2,i_3,i_4\} \subset \{1,2,\dots,n\}$, the points
\begin{equation}
    z + e_{i_1} + e_{i_2}, \quad 
    z + e_{i_2} + e_{i_3}, \quad 
    z + e_{i_3} + e_{i_4}, \quad 
    z + e_{i_1} + e_{i_4}, \quad 
    z + e_{i_1} + e_{i_3}   \label{eq:ws}
\end{equation}
belong to $L$. Then there exists a BTB graph $\tilde{\pb}$ move-equivalent to $\pb$ containing white vertices 
$w_{12}$, $w_{23}$, $w_{34}$, $w_{14}$, $w_{13}$ and black vertices 
$b_1, b_3$ such that
\[
  \ban(b_1) = z + e_{i_1}, 
    \qquad 
    \ban(b_3) = z + e_{i_3}, \quad  \ban(w_{rs}) = z + e_{i_r} + e_{i_s}
\]
that form a resplit configuration.
\end{lemma}

\begin{proof}
The points \eqref{eq:ws} form a weakly separated collection and therefore can be enlarged to a maximal weakly separated collection which corresponds to a BTB graph $\tilde \pb$ by Theorem~\ref{thm:independentsets}. Then 
\[
\{z+e_{i_1}+e_{i_4}, z + e_{i_1} + e_{i_4}, z + e_{i_1} + e_{i_3}\}~\text{and}~\{ z + e_{i_2} + e_{i_3}, 
    z + e_{i_3} + e_{i_4}, z + e_{i_1} + e_{i_3}\}
\]
are subsets of black cliques which give the black vertices $b_1,b_3$.
\end{proof}

The connection to TCD maps is as follows.

\begin{theorem}\label{th:tcddskp}
Let $\tcd$ be a minimal TCD with BTB graph $\pb$, and let 
$L \subset A_{n-1}^{k-1}$ denote the set of $\ban$-labels of all white vertices appearing in BTB graphs that are move-equivalent to $\pb$. Then every TCD map 
\[
\tcdmap : W \to \CP^\dimension
\]
on a \plabic{} graph $\pb$ with $n$ strands is the restriction of a dSKP lattice 
\[
\mathsf{S} : L \longrightarrow \CP^\dimension
\]
to the subset $\ban(W) \subset L$. Moreover, if $\tilde{\tcdmap}$ is a TCD map on a BTB graph $\tilde{\pb}$ 
obtained from $\tcdmap$ by a sequence of local moves, then 
$\tilde{\tcdmap}$ is the restriction of the \emph{same} dSKP lattice $\mathsf{S}$ 
to $\ban(\tilde{W})$, where $\tilde{W}$ is the set of white vertices of $\tilde{\pb}$.
\end{theorem}
\begin{proof}
By Lemma~\ref{lem:whitestar}, five vertices of an octahedron form a resplit configuration, and the remaining vertex appears upon performing the resplit.
Hence, by Proposition~\ref{prop:resplitmr}, the dSKP equation holds.
\end{proof}

\section{Projective and affine cluster structures}
\label{sec:clusterstructures}

In this section we endow TCD maps with two complementary cluster structures. The first, the projective cluster structure, is inherited from VRCs and is invariant under projective transformations. The second is new: the affine cluster structure depends on a choice of hyperplane at infinity and is invariant under affine transformations.

\subsection{Cluster algebras}

We begin with a brief reminder on cluster algebras. Cluster algebras are a class of commutative algebras that were introduced by Fomin and Zelevinsky~\cite{MR1887642}. Unlike typical presentations of algebras in terms of generators and relations, a cluster algebra is described by providing a set of generators as follows. We start with an initial finite set of generators and produce recursively new generators via a combinatorial operation called mutation. In this section, we give an informal introduction to cluster algebras in the restricted setting arising from \plabic{} graphs.,

\begin{definition}[Quivers]
    A \emph{quiver} is an oriented graph without self-loops and oriented $2$-cycles. We denote the vertices and \emph{arrows} (oriented edges) of a quiver $\quiver$ by $\quiver_0$ and $\quiver_1$ respectively. The vertex set $\quiver_0$ is partitioned into
    \[
    \quiver_0 = \quiver_0^{\text{mut}} \sqcup \quiver_0^{\text{fr}},
    \]
    where the vertices in $\quiver_0^{\text{mut}}$ (resp.~$\quiver_0^{\text{fr}}$) are called \emph{mutable} (resp.~\emph{frozen}).
\end{definition}

In the present article, we only deal with the following special class of quivers.

\begin{definition}[Planar dual-bipartite quivers]
     A \emph{planar dual-bipartite quiver} (or \emph{PDB quiver}) is a quiver $\quiver$ with a planar embedding such that around each vertex the arrows are alternatingly pointing outwards and inwards.
\end{definition}

As a consequence of the definition, the faces of a PDB quiver may be be colored in a bipartite manner, that is, adjacent faces get different colors (for this, we do not count the outer face as a face of the quiver). Equivalently, the dual graph of a PDB quiver is a bipartite graph.

\begin{figure}
	\includegraphics{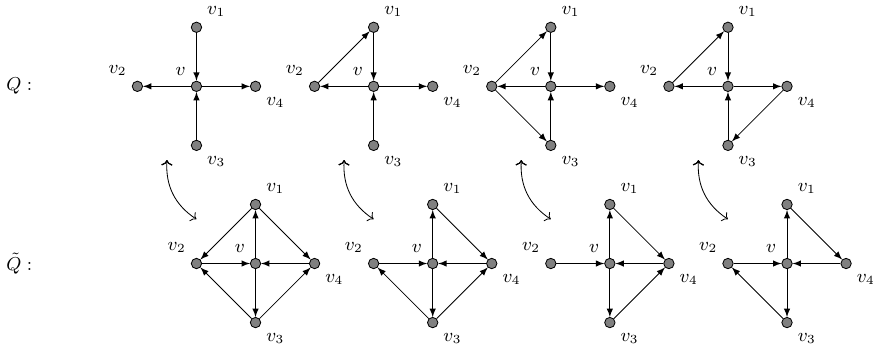}
	
	\caption{Mutation at a degree-4 vertex $v$ in a PDB quiver.}
	\label{fig:quivermut}
\end{figure}

\begin{definition}[Mutation of a PDB quiver]
    Let $\quiver$ be a PDB quiver, and let $v \in \quiver_0$ be a mutable vertex of degree $4$ with cyclically ordered neighbors $v_1, v_2, v_3, v_4$, such that the arrows incident to $v$ are
    \[
    v_1 \to v, \quad v \to v_2, \quad v_3 \to v, \quad v \to v_4.
    \]
    The \emph{mutation} of $\quiver$ at $v$ produces a new quiver $\tilde \quiver$ as follows:
    \begin{itemize}
        \item Reverse all arrows incident to $v$.
        \item Add new arrows
        \[
        v_1 \to v_2, \quad v_1 \to v_4, \quad v_3 \to v_2, \quad v_3 \to v_4.
        \]
        \item Remove any oriented $2$-cycles that are created.
    \end{itemize}
    See Figure~\ref{fig:quivermut} for an illustration.
\end{definition}

\begin{remark}
    In this paper, we restrict to mutations at degree-$4$ mutable vertices only, as these correspond to spider moves in the associated \plabic{} graph. On the other hand, in the general theory of cluster algebras, mutations are defined at all vertices of a quiver. A mutation at an arbitrary vertex of a PDB quiver may give quivers that are not PDB.
\end{remark}

It is also possible to encode the orientations of the arrows in $\quiver$ via an antisymmetric bilinear form
\begin{align*}
    \nu \colon \mathbb{Z}^{\quiver_0} \times \mathbb{Z}^{\quiver_0} &\longrightarrow \mathbb{Z}\\ (v_1, v_2) &\longmapsto \#\{\text{arrows } v_1 \to v_2\} - \#\{\text{arrows } v_2 \to v_1 \}.
\end{align*}

Then $\nu$ transforms under mutation at a vertex $v$ according to the rule:
\[
    \tilde{\nu}({v_1, v_2}) =
    \begin{cases}
        -\nu({v_1, v_2}) &\text{if $v \in \{ v_1, v_2 \}$},  \\
          \nu({v_1, v_2}) + \max(0, \nu({v_1, v})) \max(0, \nu({v, v_2})) - \max(0, \nu({v_2, v})) \max(0, \nu({v, v_1})) & \text{otherwise}.
    \end{cases}
\]

\begin{definition}[$X$-cluster variables and their mutations]\label{def:mutationx}
	The \emph{$X$-cluster variables} of a quiver $\quiver$ are a collection of variables $X = \{X_v\}_{v\in Q_0^{\text{mut}}}$ associated to the mutable vertices and taking values in $\mathbb{C}^*$. The $X$-cluster variables transform under mutation at a vertex $v$ according to the following rule: 
	\begin{equation}
\tilde X_{v'} = \begin{cases}
    X_v^{-1} & \text{if } v' = v,\\
    X_{v'} (1+X_v)^{\nu(v,v')} & \text{if } \nu(v,v')>0,\\
    X_{v'} (1+X_v^{-1})^{-\nu(v',v)} & \text{if } \nu(v',v)>0,\\
    X_{v'} & \text{otherwise.}
\end{cases}
\label{eq:mutationx}
\end{equation}
    A pair $(\quiver,X)$ is called an \emph{$X$-seed}.
\end{definition}

By starting from an initial seed and iterating mutations at every possible mutable vertex (regardless of its degree), one can generate all $X$-seeds of the cluster algebra, which may be infinite in number. Thus, a single initial $X$-seed $(\quiver, X)$ determines the entire structure, and we refer to this data as a \emph{cluster structure}.

\begin{remark}
    The $X$-cluster variables define coordinates on an algebraic torus
    \[
        \mathcal{X}_\quiver := (\mathbb{C}^*)^{\quiver_0^{\text{mut}}},
    \]
    which comes equipped with a canonical Poisson structure. The mutation formula~\eqref{eq:mutationx} defines Poisson birational transition maps between the tori $\mathcal{X}_\quiver$ and $\mathcal{X}_{\tilde\quiver}$. These maps can be used to glue together the tori corresponding to all quivers obtained from an initial quiver by mutation, producing a space $\mathcal{X}$ known as a \emph{cluster Poisson variety}~\cite{fghighertm}.  
    While the cluster variety and its Poisson structure on $\mathcal{X}$ plays a crucial role in the connection between TCD maps and cluster integrable systems studied in~\cite{agrcrdyn}, it will not play a role in the present paper.
\end{remark}

\subsection{The projective cluster structure of a TCD map}
\label{sec:projcluster}

A cluster structure for VRCs was introduced in \cite{vrc}; it induces the projective cluster structure of a TCD map.

\begin{definition}[Projective cluster structure of a TCD map]\label{def:projclusterstructure}
    Let $\tcdmap$ be a TCD map on a \plabic{} graph $\pb$. Its \emph{projective cluster structure} $\pro(\tcdmap) = (\projectivequiver, X)$ is defined as follows. The \emph{projective quiver} $\projectivequiver$ is obtained as the oriented dual of $\pb$ as follows (see Figure~\ref{fig:quivers}~(left) for an illustration):
    \begin{itemize}
        \item Place a mutable vertex of $\projectivequiver$ in each internal face of $\pb$, and a frozen vertex in each boundary face.
        \item For each edge in $\pb$, draw a dual arrow between the adjacent faces in the quiver.
        \item Orient this arrow so that it goes counterclockwise around black vertices and clockwise around white vertices.
        \item Delete any resulting oriented $2$-cycles.
        \end{itemize}
    Let $(\vrc,\relation)$ be a VRC of $\tcdmap$. For $f \in F$ a face of $G$ with boundary vertices $b_1,w_1,\dots,b_m, w_m$ in counterclockwise cyclic order, the projective cluster variable $X_f$ is defined as
\[
		X_f := (-1)^{m+1}\prod_{i=1}^m \frac{ \mu(b_i w_i)}{\mu(b_{i} w_{i-1})}.  \qedhere
\]
\end{definition}

Because the projective cluster variables are defined as alternating ratios, the variables $X_f$ associated to interior faces $f$ are invariant under gauge transformations. In other words, the cluster variables of the projective cluster structure are independent of the choice of VRC $(\vrc, \relation)$. Once a cluster structure has been defined, the next step is to understand how mutations arise. 

\begin{figure}
	\centering
	\raisebox{-0.5\height}{\includegraphics[scale=1]{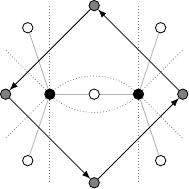}}
	\hspace{2mm}$\leftrightarrow$\hspace{1mm}
	\raisebox{-0.5\height}{\includegraphics[scale=1]{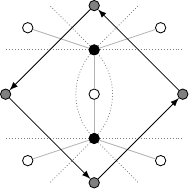}}
	\hspace{12mm}	
	\raisebox{-0.5\height}{\includegraphics[scale=1]{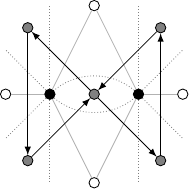}}
	\hspace{2mm}$\leftrightarrow$\hspace{1mm}
	\raisebox{-0.5\height}{\includegraphics[scale=1]{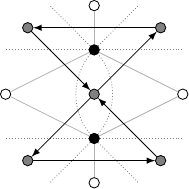}}

	\caption{The projective quiver (gray vertices and arrows) before and after a \resplitmove{} (left) as well as before and after a spider move (right).}\label{fig:twotwoquiver}
\end{figure}

\begin{lemma}\label{lem:moves_proj_cluster_structure}
	A spider move in $\pb$ corresponds to a mutation in the projective cluster structure, while a \resplitmove{} in $\pb$ leaves the projective cluster variables invariant (see~Figure~\ref{fig:twotwoquiver}).
\end{lemma}

\begin{proof}
	Examining the edge weights involved in the spider move (see Figure~\ref{fig:tcd_move_spider}), we observe that the cluster variable $X$ at the face where the mutation occurs satisfies
	\begin{align*}
		X = -\frac{bd}{ac} = \tilde{X}^{-1},
	\end{align*}
	 matching the mutation rule for $X$-cluster variables.

	For the cluster variables of adjacent faces, consider for instance the face incident to $w_1$ and $w_2$, whose $X$-variable we denote by $X_{12}$. We compute:
	\begin{align*}
		\frac{\widetilde{X}_{12}}{X_{12}} = \frac{ab^{-1} - dc^{-1}}{ab^{-1}} = 1 - \frac{bd}{ac} = 1 + X,
	\end{align*}
	which again matches the expected mutation formula. Similar computations hold for the other adjacent $X$-variables $\widetilde{X}_{23}$, $\widetilde{X}_{34}$, and $\widetilde{X}_{14}$.
	Therefore, the spider move acts as a mutation on the projective cluster structure.

	For the \resplitmove{} (Figure~\ref{fig:tcd_move_resplit}~(left)), we can verify immediately by inspection that the $X$-cluster variables remain unchanged for all four adjacent faces. 
\end{proof}

Since TCD maps are equivalent to VRCs modulo gauge, and since the $X$-cluster variables are gauge-invariant, it is natural to look for a direct definition of $X_f$ in terms of the TCD map. We begin by relating the edge weights to oriented length-ratios in affine gauge.

\begin{lemma}\label{lem:edgeratios}
	Assume we have a trivalent black vertex $b$ in a \plabic{} graph $\pb$ with neighbors $w_1,w_2,w_3$. Let $\edgeweight({b w_1}),\edgeweight({bw_2}),\edgeweight({b w_3})$ be the corresponding edge weights in affine gauge. Then, the oriented length-ratio (recall Definition~\ref{def:orientedlengthratio}) is
	\begin{equation}\label{eq:olr}
        \lambda(\tcdmap(w_2), \tcdmap(w_3), \tcdmap(w_1)) = \frac{\vrc(w_2)-\vrc(w_1)}{\vrc(w_3)-\vrc(w_1)} = -\frac{\edgeweight({bw_3})}{\edgeweight({bw_2})}.\qedhere
	\end{equation}
\end{lemma}
\begin{proof}
    Let us denote $\edgeweight({bw_i})$ by $\edgeweight_i$ for $i = 1,2,3$. Since we are in the affine gauge of a projective line, $\vrc(w_1),\vrc(w_2),\vrc(w_3)$ and $\edgeweight_1,\edgeweight_2,\edgeweight_3$ are nonzero complex numbers satisfying the following two equations:
	\begin{align*}
		\edgeweight_1 \vrc(w_1) + \edgeweight_2 \vrc(w_2) + \edgeweight_3 \vrc(w_3) &= 0,\\
		\edgeweight_1 + \edgeweight_2 + \edgeweight_3  &= 0.
	\end{align*}
	Substituting into \eqref{eq:olr}, we obtain:
    \[
        \frac{\vrc(w_2) - \vrc(w_1)}{\vrc(w_3) - \vrc(w_1)}
        = \frac{\frac{\edgeweight_1 \vrc(w_1) + \edgeweight_3 \vrc(w_3)}{\edgeweight_1 + \edgeweight_3} - \vrc(w_1)}{\vrc(w_3) - \vrc(w_1)}
        = \frac{\frac{-\edgeweight_3 \vrc(w_1) + \edgeweight_3 \vrc(w_3)}{\edgeweight_1 + \edgeweight_3}}{\vrc(w_3) - \vrc(w_1)}
        = -\frac{\edgeweight_3}{\edgeweight_2}.\qedhere
    \]
\end{proof}

Moreover, since the $X$-variables are products of certain oriented length-ratios, it turns out they may be expressed as multi-ratios.

\begin{proposition}\label{prop:projclusterviadistances}
	Let $\tcdmap$ be a TCD map. Consider a face $f \in F$ with boundary vertices $b_1,w_1,\dots,b_m, w_m$ in counterclockwise cyclic order. For each black vertex $b_i$, let $v_i$ be the third neighbor besides $w_{i-1}$ and $w_i$. Then,
    \begin{align*}
		X_f &= (-1)^{m+1} \mr(\tcdmap(w_{1}), \tcdmap(v_{2}), \tcdmap(w_{2}), \tcdmap(v_{3}), \dots, \tcdmap(w_{m}), \tcdmap(v_{1})) \\
		    &= - \prod_{i=1}^{m} \lambda(\tcdmap(w_{i-1}),\tcdmap(w_{{i}}), \tcdmap(v_{{i}}) ).\qedhere
	\end{align*}

\end{proposition}
\begin{proof}
    Direct consequence of Lemma~\ref{lem:edgeratios}. 
\end{proof}

Applying Lemma~\ref{lem:projinvariants} to Proposition~\ref{prop:projclusterviadistances} we get the following observation, justifying the terminology \emph{projective} cluster structure.
\begin{corollary}
    The projective cluster variables are invariant under projective transformations.
\end{corollary}

\begin{remark}
Another important observation is that the quiver $\projectivequiver$ determines $\pb$ up to \resplitmove{}s. Each counterclockwise face of the quiver contains exactly one white vertex. 
The restriction of $\pb$ to a clockwise face of $\projectivequiver$ is a tree, with all possible trees related by \resplitmove{}s. 
Consequently, each clockwise face $f$ of $\projectivequiver$ of degree $\deg(f)$ contains $\deg(f) - 2$ black vertices and $\deg(f) - 3$ white vertices of $\pb$. 
These clockwise faces encode the geometric condition that the points of $\tcdmap$ corresponding to the adjacent white faces lie in a $(\deg(f) - 2)$-dimensional projective subspace.
\end{remark}

We now describe the space of TCD maps with given projective cluster variables. Recall that $F^\circ$ denotes the set of interior faces of the \plabic{} graph $G$.
The following result was proved in \cite{vrc} (in different terminology).

\begin{theorem}[{\cite[Theorem~7.8]{vrc}}] \label{thm:X_exists}
    Let $\tcd$ be a TCD of maximal rank. For any set of projective cluster variables $X:F^\circ\to\C^*$, there exists a unique TCD map (up to $\operatorname{PGL}(\maxrank(\tcd)+1)$) that realizes $X$.
\end{theorem}

Let us emphasize that, once $X$ is fixed, there is no freedom to prescribe boundary points (up to $\operatorname{PGL}(\maxrank(\tcd)+1)$). More generally, we show the following.

\begin{proposition}
    Let $\tcd$ be a TCD with \plabic{} graph $\pb$ and let $1 \leq \dimension \leq \maxrank(\tcd)$. For every assignment of projective cluster variables $X:F^\circ \rightarrow \C^*$, there exists a TCD map $\tcdmap \in \maps(\tcd,\dimension)$ realizing $X$. Moreover, for fixed $X$, the set of TCD maps of rank $\dimension$ modulo $\operatorname{PGL}(\dimension+1)$ realizing $X$ is a variety of dimension 
    \[
        (\maxrank(\tcd)-\dimension)(\dimension+1). \qedhere
    \]
\end{proposition}

\begin{proof}
    By Theorem~\ref{thm:X_exists}, there is a unique maximal-rank map 
    $\hat\tcdmap$ realizing $X$ modulo $\operatorname{PGL}(\maxrank(\tcd)+1)$. Let $(\hat\vrc,\hat\relation)$ be a VRC lift of $\hat\tcdmap$. By a projective transformation and a gauge transformation, we may assume that the vectors 
    \[
        \hat \vrc(w_i^\partial), \qquad w_i^\partial \in W^{\min},
    \]
    are the standard basis vectors for $\C^{\maxrank(\tcd)+1}$.

    Since central projections preserve multi-ratios (Lemma~\ref{lem:projinvariants}), any generic central projection of $\hat\tcdmap$ is a TCD map that realizes $X$. 

    Consider the surjective rational map 
    \begin{align*}
        \Phi:\operatorname{Gr}(\maxrank(\tcd)-\dimension,\maxrank(\tcd)+1) &\dashrightarrow  \{\text{rank-}\dimension \text{ TCD maps realizing }X\}/\operatorname{PGL}(\dimension+1)\\
        \subspacea &\longmapsto \pi_{\subspacea}\circ \hat\tcdmap,
    \end{align*}
    defined on the Zariski-open subset of $\hat \tcdmap$-admissible $\subspacea$. In homogeneous coordinates, by \eqref{eq:central_projection_affine}, $\Phi$ is given by the projectivization of 
    \[
\tilde \subspacea \longmapsto \tilde \pi_{\tilde \subspacea} \circ \hat \vrc.
    \]
    We now show that the map $\Phi$ is birational. Given $\tcdmap \in \maps(\tcd,\dimension)$ realizing $X$, by Proposition~\ref{prop:central_projection_lift_exists}, there is a $\subspacea$ and a VRC lift $(\vrc,\relation)$ of $\tcdmap$ such that $\vrc = \tilde \pi_{\tilde\subspacea} \circ \hat \vrc$. Moreover, the subspace $\subspacea$ can be recovered from the map 
    \begin{align*}
        \C^{\maxrank(\tcd)+1} &\longrightarrow \C^{\dimension+1}, \\
        \hat \vrc(w_i^\partial) &\longmapsto \vrc(w_i^\partial),
    \end{align*}
    as the projectivization of its kernel. Therefore $\Phi$ is birational, and thus
    \begin{align*}
       \dim \{\text{rank-}\dimension \text{ TCD maps realizing }X\}/\operatorname{PGL}(\dimension+1) &= \dim \operatorname{Gr}(\maxrank(\tcd)-d,\maxrank(\tcd)+1)\\&=(\maxrank(\tcd)-d)(d+1).\qedhere
    \end{align*}
\end{proof}

\subsection{The affine cluster structure of a TCD map}\label{sec:affcluster}

\begin{figure}
	\centering
	\raisebox{-0.5\height}{\includegraphics[scale=1.1]{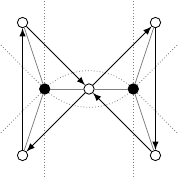}}    
	\hspace{2mm}$\leftrightarrow$\hspace{1mm}
	\raisebox{-0.5\height}{\includegraphics[scale=1.1]{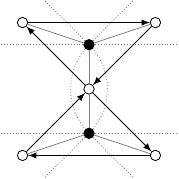}}
	\hspace{10mm}	
	\raisebox{-0.5\height}{\includegraphics[scale=1.1]{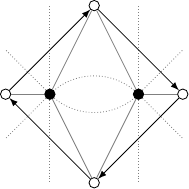}}
	\hspace{2mm}$\leftrightarrow$\hspace{1mm}
	\raisebox{-0.5\height}{\includegraphics[scale=1.1]{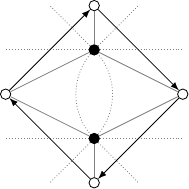}}

	\caption{The transformation of the quiver of the affine cluster structure for the \resplitmove{} (left) and the spider move (right).}\label{fig:twotwoaffinequiver}
\end{figure}

We now introduce a second type of cluster structure for TCD maps, which is new.

\begin{definition}[Affine cluster structure of a TCD map] \label{def:affinecluster}
    Let $\hyperplane$ the hyperplane at infinity for an affine chart. Let $\tcdmap$ be a TCD map, and let $\mu$ be the edge weights of an associated VRC in affine gauge. Assume that $\hyperplane$ is generic with respect to $\tcdmap$. The \emph{affine cluster structure} $\aff_{\hyperplane}(\tcdmap) = (\affinequiver, Y)$ is defined as follows. The \emph{affine quiver} $\affinequiver$ is constructed via the following steps (see Figure~\ref{fig:quivers}~(center) for an illustration):
    \begin{enumerate}
        \item Place a mutable vertex of $\quiver^{\mathrm{aff}}$ at each internal white vertex of $\pb$, and a frozen vertex at each boundary white vertex.
        \item For each trivalent black vertex with neighbors $w_1,w_2,w_3$ in clockwise cyclic order, draw arrows  
        \[
        w_1 \to w_2, \quad w_2 \to w_3, \quad w_3 \to w_1.
        \]
        \item For each $i$, draw an arrow $w_{i+1}^\partial \to w_i^\partial$.
        \item Delete any resulting oriented $2$-cycles.
    \end{enumerate}

	Let $w$ be an internal white vertex of $\pb$, and suppose $b_1, \dots, b_m$ are the black vertices adjacent to $w$, ordered in counterclockwise order around $w$. For each $b_i$, let $w_i$ and $w'_i$ be the other two white neighbors of $b_i$, such that the counterclockwise cyclic order around $b_i$ is $w, w_i, w'_i$. Then the \emph{affine cluster variable} associated to $w$ is defined as the alternating ratio:
	\begin{align*}
		Y_w := (-1)^{m+1} \prod_{i=1}^m \frac{\mu(b_i w'_i)}{\mu(b_i w_i)}.
	\end{align*}
\end{definition}

\begin{figure}
	\centering
	\raisebox{-0.5\height}{\includegraphics[scale=1.]{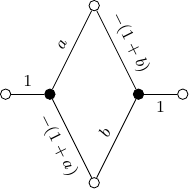}}    
	\hspace{3mm}$\leftrightarrow$\hspace{0mm}		
	\raisebox{-0.5\height}{\includegraphics[scale=1.]{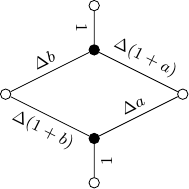}}    
	\hspace{10mm}
	\raisebox{-0.5\height}{\includegraphics[scale=1.]{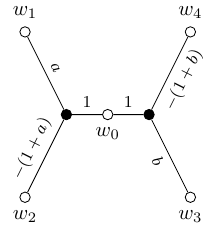}}    
	\hspace{0mm}$\leftrightarrow$\hspace{0mm}
	\raisebox{-0.5\height}{\includegraphics[scale=1.]{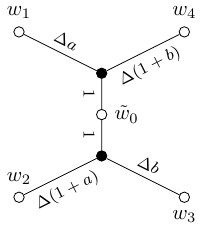}}    

	\caption{Change of edge weights in affine gauge, with $\Delta=-(1+a+b)^{-1}$.}\label{fig:affineweights}
\end{figure}

Next, we study the effect of moves on the affine cluster structure.

\begin{lemma} \label{lem:affine_cluster_moves}
	The \resplitmove{} corresponds to a mutation in the affine cluster structure, while the spider move leaves the affine cluster structure invariant (see~Figure~\ref{fig:twotwoaffinequiver}).
\end{lemma}

\begin{proof}
	Figure~\ref{fig:affineweights} shows the change of edge weights under the two moves. The affine cluster variables are clearly unchanged under the spider move for each of the four involved white vertices.

	For the \resplitmove{} at a white vertex $w_0$, we compute:
	\begin{align*}
		Y_{w_0} = - (1 + a^{-1})(1 + b^{-1}) = \tilde{Y}_{\tilde w_0}^{-1},
	\end{align*}
	which matches the mutation rule.

	For the resplit at $w_2$, we have:
	\begin{align*}
		\frac{\tilde{Y}_{w_2}}{Y_{w_2}} = \frac{1}{a\Delta b} = - (a^{-1} + b^{-1} + a^{-1}b^{-1}) = 1 + Y_{w_0},
	\end{align*}
    which also matches the mutation rule. The calculations for the other adjacent white vertices are similar.
\end{proof}

Recall that we were able to write the projective cluster variables as products of oriented length-ratios. Similarly, we show that the affine cluster variables can be written in terms of (other) products of oriented length-ratios.

\begin{definition}[Star-ratio]\label{def:starratio}
	Let $\point,\point_1,\point'_1,\point_2,\point'_2,\dots, \point_m,\point'_m$ be $2m+1$ points in $\CP^\dimension$ such that $\point$ lies on each line ${\point_i \point'_i}$. Choose an affine chart $\C^{\dimension} \subset \CP^\dimension$ containing all the points and define the \emph{star-ratio} of these points as 
\[
        \sr(\point;\point_1,\point'_1,\point_2,\point'_2\dots,\point_m,\point'_m) := \frac{\prod_{i=1}^m (\point_i-\point)}{\prod_{i=1}^m(\point-\point'_i)}
        = (-1)^m \prod_{i=1}^m \lambda(\point_i,\point'_i,\point).  \qedhere
\]  
\end{definition}
Since star-ratios are products of oriented lengths-ratios, star-ratios are invariant under affine transformations (but not general projective transformations). In particular, this means that the affine cluster variables of a TCD map in two different affine charts but with the same hyperplane at infinity are equal.
Using the notion of star-ratios, we are able to give a geometric definition of the affine cluster variables.
\begin{proposition}\label{lem:starratioviadistances}
	Let $\tcdmap$ be a TCD map in affine gauge with respect to a generic hyperplane $\hyperplane$. Then the affine cluster variable at a white vertex $w$ is given in terms of a star-ratio in this affine chart as
	\begin{align*}
		Y_w &= (-1)^{m+1} \sr(\tcdmap(w);\tcdmap(w_1),\tcdmap(w_1'),\tcdmap(w_2),\tcdmap(w_2'),\ldots,\tcdmap(w_m),\tcdmap(w_m'))\\
        &=-\prod_{i=1}^{m} \lambda(\tcdmap(w_i),\tcdmap(w_i'),\tcdmap(w)),
	\end{align*}
    where $w_1,w'_1,\dots,w_m,w'_m$ are as in Definition~\ref{def:affinecluster}. 
\end{proposition}
\begin{proof}
	This follows directly from Lemma~\ref{lem:edgeratios}.
\end{proof}

Since star-ratios are invariant under affine transformations, we immediately obtain the following corollary.

\begin{corollary}
    The affine cluster variables are invariant under affine transformations.
\end{corollary}
Thus the terminology of affine cluster structure is justified.

\begin{remark}
    As in the case of the projective cluster structure, the quiver determines the \plabic{} graph, but in this case, up to spider moves. Since the white vertices of $\pb$ are the same as the vertices of $\affinequiver$, we can place the points of the TCD map at the vertices of $\affinequiver$. Then each clockwise face of the quiver corresponds to a common line shared by the $\deg(f) - 2$ black vertices of $\pb$ contained in that face. The white vertices corresponding to the vertices of $\affinequiver$ adjacent to a clockwise face of the quiver all lie on this common line.
\end{remark}

\subsection{Affine cluster structure and t-embeddings}\label{sec:tembeddings}

In this section, we describe the cluster structure of t-embeddings and show that it is an example of an affine cluster structure. Let $\quiver$ be a PDB quiver and consider a map 
\[
\mathsf{t} : \quiver_0 \longrightarrow \C 
\]
where $\C \subset \CP^1$ denotes the standard affine chart. For $v \in \quiver_0$, define
\begin{equation}
       Y_v := -\prod_{i=1}^{m} \lambda(\mathsf{t}(v_i),\mathsf{t}(v_{i}'),\mathsf{t}(v)), \label{eq:Y_var_t_embedding} 
\end{equation}
where $v_1, v_1',\dots,v_m,v_m'$ are the neighbors of $v$ in counterclockwise order such that $v_i \to v$ and $v \to v_i'$ are the edges incident to $v$. We call $\mathsf{t}$ a \emph{t-realization} if $Y_v>0$ for all $v \in \quiver_0$ and a \emph{t-embedding} if, additionally, the straight-line drawing of $\quiver$ determined by $\mathsf{t}$ is an embedding. Note that adjacent points of a t-realization cannot coincide since otherwise some $Y_v$ would be 0 or $\infty$. t-realizations are called \emph{planar conical nets} in~\cite{amiquel,muellerconical} and \emph{Coulomb gauges} in~\cite{kenyonlam}, though the term {t-realization/embedding}, introduced in~\cite{clrtembeddings}, has become the most widely used.

It was shown in~\cite{amiquel,kenyonlam} that the $Y_v$ are the cluster variables of a cluster structure with quiver $\quiver$. We now show that it is the affine cluster structure of an associated TCD map. Construct a BTB graph $\pb$ from $\quiver$ as follows: 
\begin{itemize}
    \item Place a white vertex $w_v$ at each vertex $v \in \quiver_0$ and let $Y_{w_v}:=Y_v$.
    \item Triangulate each clockwise face of $\quiver$ and replace each triangle as follows:
    
\begin{center}
\begin{tikzpicture}[scale=1,baseline={(current bounding box.center)}]

\begin{scope}[shift={(0,0)}]
  \coordinate (a) at (90:1);
  \coordinate (b) at (210:1);
  \coordinate (c) at (330:1);
  \draw[gray, dashed] (a)--(b)--(c)--cycle;
\end{scope}
\node at (2.5,0) {$\longmapsto$};

\begin{scope}[shift={(5,0)}]
  \coordinate[wvert] (v1) at (90:1);
  \coordinate[wvert] (v2) at (210:1);
  \coordinate[wvert] (v3) at (330:1);
  \draw[gray, dashed] (v1)--(v2)--(v3)--(v1);
  \coordinate[bvert] (o) at (0,0);
  \draw (o)--(v1) (o)--(v2) (o)--(v3);
\end{scope}
\end{tikzpicture}.
\end{center}
  
\end{itemize}
Given a t-realization $\mathsf{t}$, define a TCD map $\tcdmap: W \longrightarrow \CP^1$ by
\[
    \tcdmap(w_v) := \mathsf{t}(v).
\]
Since $\tcdmap$ maps to $\CP^1$, the incidence conditions are trivially satisfied. 
\begin{proposition}
    The cluster structure $(\quiver,Y)$ is the affine cluster structure of $\tcdmap$.
\end{proposition}
\begin{proof}
By the construction of $\pb$ from $\quiver$, the quiver $\quiver$ is precisely the affine quiver associated with $\pb$. Moreover, by Proposition~\ref{lem:starratioviadistances}, the expression~\eqref{eq:Y_var_t_embedding} is the corresponding affine cluster variable. 
\end{proof}

\begin{remark}
Since all points lie in $\CP^1$, there is a curious symmetry: one could equally obtain a TCD map with the same affine cluster variables by triangulating the counterclockwise faces of $\quiver$ to construct $\pb$.
\end{remark}

\section{Sections of TCD maps}
\label{sec:sections}

In this section, we investigate what happens when we intersect a TCD map with a hyperplane.  Each line of the TCD map becomes a point and these points satisfy new incidence relations giving rise to a new TCD map. Before defining this operation, we introduce a genericity condition.

\begin{definition}[$1$-generic TCD maps]
    A TCD map $\tcdmap$ is \emph{$1$-generic} if:
    \begin{itemize}
        \item For every face $f$ of $\pb$ and three white vertices $w_1,w_2,w_3$ on the boundary of $f$, the points $\tcdmap(w_1)$, $\tcdmap(w_2)$, $\tcdmap(w_3)$ span a plane.
        \item For every consecutive pair of boundary white vertices $w_i^\partial, w_{i+1}^\partial$ that are distinct (not glued in the cactus), the points $\tcdmap(w_i^\partial), \tcdmap(w_{i+1}^\partial)$ are also distinct. \qedhere
    \end{itemize}
\end{definition}

\begin{figure}
	\centering
	\hspace{-10mm}
	\raisebox{-0.5\height}{\includegraphics[scale=0.9]{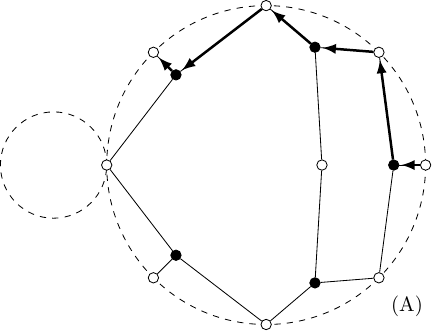}}
	\hspace{10mm}
	\raisebox{-0.5\height}{\includegraphics[scale=0.9]{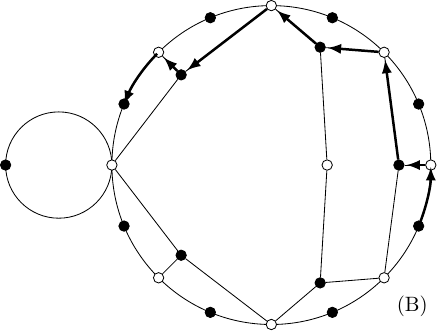}}
    \vspace{-1mm}
	\\
	\hspace{15mm}
	\raisebox{-0.5\height}{\includegraphics[scale=0.9]{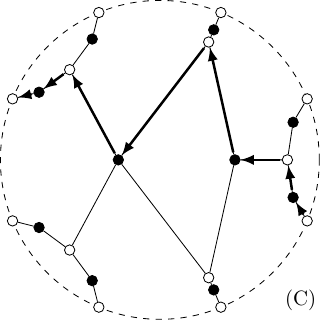}}    
	\hspace{18mm}
	\raisebox{-0.5\height}{\includegraphics[scale=0.9]{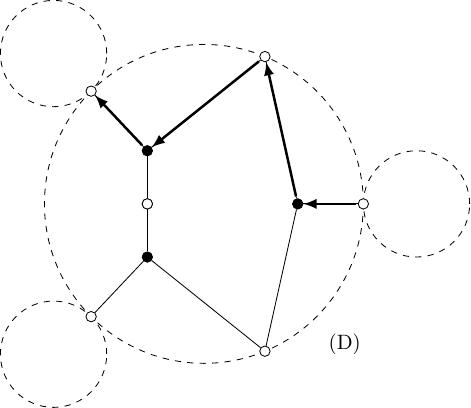}}    
	\caption{(A) A \Plabic{} graph $\pb$; (B) the graph $\pb^{\mathrm{int}}$; (C) the star graph $\stargraph$; (D) the section $\sigma(\pb)$. We also illustrate the correspondence of zig-zag paths by tracking a single zig-zag path across all four graphs. 
    }
	\label{fig:four-panels}
\end{figure}

Moreover, to define sections of a TCD map and study their properties, it is convenient to introduce the following auxiliary graph.

\begin{definition}[Star graph] \label{def:star_graph}
     Let $\tcd$ be a TCD with \plabic{} graph $\pb$. We construct a bipartite graph $\stargraph \subset \disk$, called the \emph{star graph} of $\pb$, as follows:
    \begin{enumerate}
		\item Place a boundary black vertex $b_{i,i+1}^{\partial}$ between boundary white vertices $w_i^\partial, w_{i+1}^\partial$ and draw edges $b_{i,i+1}^\partial w_i^\partial$ and $b_{i,i+1}^\partial w_{i+1}^\partial$ along the boundary of the cactus to get a graph $\pbint$.
		\item Place a white vertex $\starw_b$ at each black vertex $b$ of $\pbint$.
		\item Place a black vertex $\starb_f$ inside each face of $\pbint$.
		\item Draw an edge $\starb_f \starw_b$ if $b$ and $f$ are incident in $\pbint$.
        \item In the special case where $w_i^\partial,w_{i+1}^\partial$ are glued, the boundary face $f_{i,i+1}^\partial$ between them is a disk of the cactus containing no other vertices of $\pb$, and at the end of the previous step, we obtain a black leaf $\starb_{f_{i,i+1}^\partial}$ connected to $\starw_{b_{i,i+1}^\partial}$. Remove this disk of the cactus (including $\starb_{f_{i,i+1}^\partial}, \starw_{b_{i,i+1}^\partial}$ and the edge between them).
        \item Take the preimage of the resulting graph under the quotient map $\disk \rightarrow \disk/\ncp$. As the graph is disjoint from the boundary vertices $\{w_1^\partial, \dots, w_n^\partial\}$, the preimage is simply the same graph embedded in the original disk~$\disk$. \qedhere
	\end{enumerate}
\end{definition}

Figure~\ref{fig:four-panels} shows an example of the construction of the star graph. The reason for the name is that in Steps~(3)~and~(4), we place a ``star" inside each face.

\begin{proposition}\label{lem:star_graph_properties}
    The star graph $\stargraph$ satisfies the following properties:
    \begin{enumerate}
        \item The white vertices of $\pb$ correspond bijectively to the faces of $\stargraph$, and boundary white vertices correspond to boundary faces.
        \item There is a natural bijection between edges of $\pb$ and internal edges of $\stargraph$. If $bw$ is an edge of $\pb$ (with $b$ black and $w$ white) and $f$ is the face of $\pb$ opposite $w$ across $bw$, then $bw$ corresponds to the edge $\hat b_f \starw_b$ of $\stargraph$.
        \item Additionally, $\stargraph$ has $n$ boundary edges of the form $\starb_{f_{i,i+1}^\partial}\starw_{b_{i,i+1}^\partial}$, where $f_{i,i+1}^\partial$ is the boundary face of $\pb$ between $w_i^\partial$ and $w_{i+1}^\partial$. \qedhere
    \end{enumerate}
\end{proposition}
\begin{proof}
    All three statements are immediate from the construction of $\stargraph$.
\end{proof}

To define a section of a TCD map, we build upon the combinatorial construction of the star graph and add geometric data to it.

\begin{definition}[Sections of TCD maps]\label{def:section}
	Let $\tcdmap: W(\pb) \to \CP^\dimension$ with $\dimension \geq 2$ be a $1$-generic TCD map on a TCD $\tcd$ with \plabic{} graph $\pb$, and let $\hyperplane \subset \CP^\dimension$ be a hyperplane in generic position with respect to $\tcdmap$.  
	A \emph{section} of $\tcdmap$ with respect to $\hyperplane$ is a TCD map
	\[
	   \sigma_\hyperplane(\tcdmap): W(\sigma(\pb)) \longrightarrow \hyperplane \cong \CP^{\dimension-1},
	\]
	defined on a new \plabic{} graph $\sigma(\pb)$ (equivalently, TCD $\sigma(\tcd)$) constructed as follows:
	\begin{enumerate}
		\item Split every black vertex of degree greater than $3$ in $\stargraph$, replacing it with a bipartite tree with trivalent black vertices.
		\item Contract all black vertices of degree~$2$, possibly resulting in a new cactus due to identifications along the boundary, to obtain the graph $\sigma(\pb)$.
		\item By Proposition~\ref{lem:star_graph_properties} and because contractions and splits induce bijections on faces, the faces of $\sigma(\pb)$ correspond bijectively to white vertices of $\pb$; we denote the face corresponding to $w$ by $\starf_w$. For a white vertex $\starw$ of $\sigma(\pb)$, set
		\[
	       	\projectiveline(\starw) := \Span\bigl\{\tcdmap(w): \starf_w \text{ is incident to } \starw\bigr\}.
		\]
		By $1$-genericity of $\tcdmap$, each $\projectiveline(\starw)$ is a line. We observe that, if $\starw_b$ is an interior white vertex of $\sigma(\pb)$ coming from a black vertex $b$ of $\pb$, then $\projectiveline(\starw_b)$ is simply the line $\projectiveline_b$ living at $b$. Define
		\[
		      \sigma_\hyperplane(\tcdmap)(\starw) := \projectiveline(\starw) \cap \hyperplane.
		\]
		Since $\hyperplane$ is generic with respect to $\tcdmap$, this intersection is a point of $\hyperplane$. \qedhere
	\end{enumerate}
\end{definition}

Note that the BTB graph $\sigma(\pb)$ does not depend on the choice of $\hyperplane$. Figure~\ref{fig:four-panels}~(D) shows the graph $\sigma(\pb)$ for the \plabic{} graph in Figure~\ref{fig:four-panels}~(A).

\begin{proposition}
\label{prop:sectiontcdmap}
Every section $\sigma_\hyperplane(\tcdmap)$ is a TCD map.
\end{proposition}
\begin{proof}
    We need to check that around every black vertex $\starb$ of $\sigma(\pb)$, the points of $\sigma_\hyperplane(\tcdmap)$ are distinct and lie on a line. The black vertices of $\sigma_\hyperplane(\tcdmap)$ are trivalent vertices of a tree in Step~(1) of Definition~\ref{def:section} inside a face $f$ of $\pb$. Let $\starf_{w_1},\starf_{w_2},\starf_{w_3}$ be the three faces incident to $\starb$, where $w_1,w_2,w_3$ are three white vertices incident to $f$. Let $\starw_{12}, \starw_{23}, \starw_{31}$ be the three white neighbors of $\starb$ so that $\starf_{w_i}$ is between $\starw_{i-1,i}$ and $\starw_{i,i+1}$.  
    
    The three points $\tcdmap (w_1),\tcdmap(w_2),\tcdmap(w_3)$ span a plane $\mathsf{\Pi}$ by $1$-genericity. Therefore the three lines 
    \[
    \projectiveline(\starw_{i,i+1}):= \tcdmap(w_i) \tcdmap(w_{i+1}), \quad i=1,2,3,
    \]
    also span $\mathsf{\Pi}$. Since $\hyperplane$ is generic with respect to $\tcdmap$, $\mathsf{\Pi}$ is not contained in $\hyperplane$ and therefore intersects it in a line $\projectiveline(\starb)$. By construction, $\projectiveline(\starb)$ contains the three points $\sigma_\hyperplane(\tcdmap)(\starw_{12}), \sigma_\hyperplane(\tcdmap)(\starw_{23}), \sigma_\hyperplane(\tcdmap)(\starw_{31})$. If two of these points are the same, say $\sigma_\hyperplane(\tcdmap)(\starw_{i,i+1}) =\sigma_\hyperplane(\tcdmap)(\starw_{j,j+1})$ for some $i \neq j$, then $\projectiveline(\starw_{i,i+1})$ and $\projectiveline(\starw_{j,j+1})$ have two points in common and therefore are the same line which then contains $\tcdmap(w_1),\tcdmap(w_2),\tcdmap(w_3)$, which contradicts $1$-genericity.
\end{proof}

Note that different choices of trees in Step~(1) of Definition~\ref{def:section} give rise to different sections. Therefore, a TCD map does not, in general, have a unique section. However, we have the following relations.

\begin{proposition} \label{lemma:section_properties}
    Let $\tcdmap$ be a 1-generic TCD map and let $\hyperplane$ be a generic hyperplane with respect to $\tcdmap$.
    \begin{enumerate}
        \item Any two sections $\sigma_\hyperplane(\tcdmap)$ are related by \resplitmove{}s.
        \item If $\tcdmap$ and $\tcdmap'$ are related by a spider move, then the sets of sections of $\tcdmap$ and $\tcdmap'$ coincide.\qedhere
    \end{enumerate}
\end{proposition}

\begin{proof}

    (1) follows directly from the construction. 
    
    For~(2), if $f$ is a quad face of $\pb$ with vertices $b_1,w_1,b_2,w_2$ in counterclockwise order, then $\projectiveline(b_1)=\projectiveline(b_2)$; let us denote this common line by $\projectiveline$. The new black vertex $\starb_f$ has degree $2$ and is contracted in Step~(2), yielding a white vertex $\starw$. The point $\projectiveline \cap H$ is the TCD-map point associated to the white vertex $\starw$ in $\sigma(\pb)$. Since a spider move does not change the line $\projectiveline$, any section is invariant under spider moves. 
\end{proof}

\begin{remark}
    $1$-genericity is stronger than what is necessary for the existence of a single section. Instead, $1$-genericity guarantees the existence of all sections.
\end{remark}

The following proposition shows how strands behave under taking a section of a TCD map.

\begin{proposition} \label{prop:section_permutation}
    Let $\tcdmap$ be a TCD map with \plabic{} graph $\pb$. Assume that no two consecutive boundary white vertices of $\pb$ are glued. Then:
    \begin{enumerate}
        \item There is a natural bijection between the strands of $\pb$ and those of $\sigma(\pb)$.
        \item The strand permutation of $\sigma(\tcd)$ is given by:
        \[
          C_{\sigma(\tcd)}(i)=C_{\tcd}(i)-1,
        \]
        with indices taken modulo $n$.
        
        \item The \plabic{} graph $\sigma(\pb)$ is minimal.\qedhere
    \end{enumerate}
\end{proposition}
\begin{proof}

    By prepending (resp.~appending) the boundary edge $b_{i,i+1}^\partial w_i^\partial$ (resp.~$w_j^\partial b^\partial_{j-1,j}$), a zig-zag path in $\pb$ from $w_i^\partial$ to $w_{j}^\partial$ corresponds to a zig-zag path in $\pbint$ from $b^\partial_{i,i+1}$ to $b^\partial_{j-1,j}$; see Figure~\ref{fig:four-panels}~(B). Now let $f$ be a face of $\pbint$ with vertices $b_1,w_1,\dots,b_m, w_m$ in counterclockwise cyclic order. 
    Then a zig-zag segment $b_i \rightarrow w_{i+1} \rightarrow b_{i+1}$ in $\pbint$ gets replaced with the segment $\starw_{b_i} \rightarrow \starb_f \rightarrow \starw_{b_{i+1}}$ in $\stargraph$. Gluing these replacements face-by-face gives a bijection between zig-zag paths in $\pbint$ and $\stargraph$; see Figure~\ref{fig:four-panels}~(C). There is a natural bijection between zig-zag paths in $\stargraph$ and $\sigma(\pb)$ induced by splits and contractions; see Figure~\ref{fig:four-panels}~(D). Composing all these bijections gives the bijection between the strands of $\pb$ and $\sigma(\pb)$ in~(1). 
    
    Under the bijection in~(1), the zig-zag path in $\pb$ from $w_i^\partial$ to $w_j^\partial$ becomes the zig-zag path in $\sigma(\pb)$ from $b_{i,i+1}^\partial$ to $b_{j-1,j}^\partial$, giving~(2).
    
    For~(3), we use the characterization of minimality in Theorem~\ref{th:mintcdforbidden}. We have already seen that every strand in $\sigma(\pb)$ starts and ends at the boundary; there are no strands that are loops in $\sigma(\pb)$.

    The correspondence between edges in $\pb$ and $\stargraph$ is compatible with zig-zag paths in the following sense: a zig-zag path of $\pb$ uses $bw$ in a given direction (black-to-white or white-to-black) if and only if the corresponding zig-zag path of $\stargraph$ uses $\starb_f \starw_b$ in the same direction. See Figure~\ref{fig:four-panels} for an illustration. Additionally, $\stargraph$ has $n$ boundary edges of the form
    $\starb_{f_{k,k+1}^\partial}\starw_{b_{k,k+1}^\partial}$, where $f_{k,k+1}^\partial$ is the boundary face of $\pb$ between
    $w_k^\partial$ and $w_{k+1}^\partial$.

    In terms of zig-zag paths, a parallel bigon arises when two zig-zag paths traverse two distinct edges in the same order. From the correspondence of zig-zag paths and edges above, no pair of internal edges of $\stargraph$ can form a parallel bigon since it would mean that the corresponding zig-zag paths form a parallel bigon in $\pb$. Furthermore the boundary edges, that is, edges of the form $\starb_{f_{k,k+1}^\partial}\starw_{b_{k,k+1}^\partial}$ are traversed exactly once in each direction by some zig-zag path: in the white-to-black (resp.~black-to-white) direction as the first (resp.~last) edge of a zig-zag path. Thus boundary edges cannot participate in parallel bigons.
  
    Similarly, a self-intersection arises when a zig-zag path traverses an edge in both directions. This cannot happen at an internal edge of $\stargraph$ since it would correspond to a self-intersection of the corresponding zig-zag path in $\pb$. At a boundary edge, a self-intersection of a zig-zag path in $\stargraph$ would require a zig-zag path of $\pb$ with $C_{\tcd}(i)=i+1$. Indeed, the corresponding zig-zag path in $\pbint$ must use the edge from $b^\partial_{i,i+1}$ to $w^\partial_i$, which can only occur at the beginning of the path, as well as the edge from $w^\partial_{i+1}$ to $b^\partial_{i,i+1}$, which can only occur at the end of the path. If $C_\tcd(i) = i+1$, the strand $i$ does not intersect any other strands. Indeed, one can construct a minimal TCD $\tcd'$ in which strand $i$ has no crossings by first taking a minimal TCD of the remaining strands and then adding strand $i$ without introducing any intersections.  
    By Theorem~\ref{th:tcdmovesconnected}, $\tcd$ is related to $\tcd'$ by 2-2 moves, and strand $i$ cannot participate in such a move.
        
    Consequently, strand $i$ must be a part of the boundary of a clockwise region of $\tcd$ that contains both $w_i^\partial$ and $w_{i+1}^\partial$. This forces these two vertices to be 
    glued, a contradiction.
\end{proof}

\begin{remark}
Recall the map $\ban$ from Section~\ref{sec:abel}. 
Via the natural bijection between white vertices of $\pb$ and faces of $\sigma(\pb)$, the map $\ban$ is preserved, i.e.
\[
  \ban(\starf_w) = \ban(w).
\]
Moreover, $\enm{n}{2}$-subconfigurations (recall Definition~\ref{def:endpointmatching}) of $\tcdmap$ correspond to white vertices of the star graph $\stargraph$. 
Within each such $\enm{n}{2}$-configuration, the map $\ban$ takes the same value on all black vertices, and this value agrees with $\ban$ of the corresponding white vertex in $\stargraph$.
\end{remark}

\section{Affine cluster structures as sections of projective cluster structures}
\label{sec:sectioncluster}

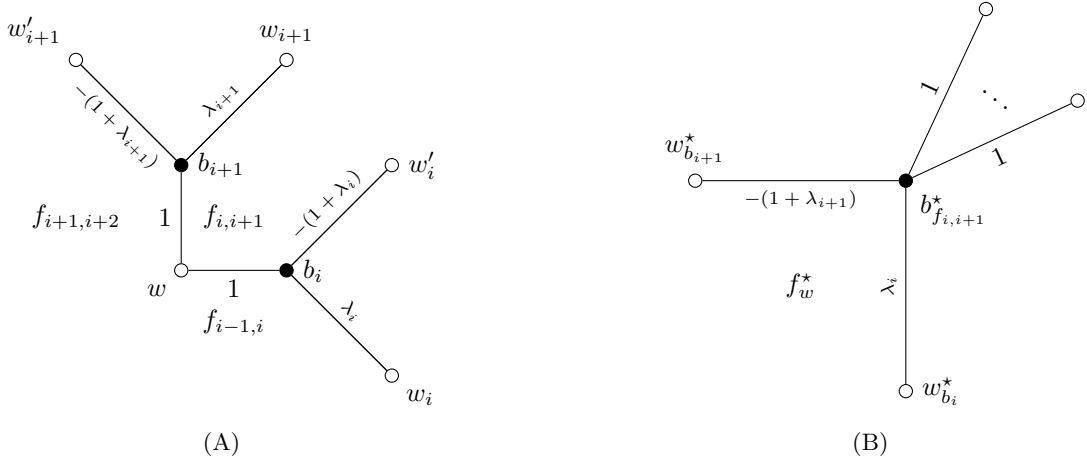
\begin{figure}
\centering

\begin{subfigure}{0.48\textwidth}
  \centering

  \begin{tikzpicture}[scale = 1.4]
    \node[wvert,label={below left}:$w$] (w) at (0,0) {};
    \node[bvert,label={right}:$b_i$] (b1) at (1,0) {};
    \node[bvert,label={right}:$b_{i+1}$] (b2) at (0,1) {};

    \node[wvert,label={below right}:$w_i$] (w1) at (2,-1) {};
    \node[wvert,label={right}:$w'_i$] (ww1) at (2,1) {};
    \node[wvert,label={above}:$w_{i+1}$] (w2) at (1,2) {};
    \node[wvert,label={above left}:$w'_{i+1}$] (ww2) at (-1,2) {};

    \node[] (f12) at (.5,.5){$f_{i,i+1}$};
    \node[] (f23) at (-1,.5){$f_{i+1,i+2}$};
    \node[] (f41) at (0.5,-.5){$f_{i-1,i}$};

    \path[font=\scriptsize]
      (b1) edge node [sloped, above] {$\lambda_i$} (w1)
           edge node [sloped, above] {$-(1+\lambda_i)$} (ww1)
      (b2) edge node [sloped, above] {$\lambda_{i+1}$} (w2)
           edge node [sloped, below] {$-(1+\lambda_{i+1})$} (ww2);

    \draw[-]
      (w)  edge node[below]{1} (b1)
      (w)  edge node[left]{1}  (b2)
      (b1) edge (w1)
      (b1) edge (ww1)
      (b2) edge (w2)
      (b2) edge (ww2);
  \end{tikzpicture}
  \caption*{(A)}
  \label{fig:section_hexagon:A}
\end{subfigure}\hfill
\begin{subfigure}{0.48\textwidth}
  \centering

  \begin{tikzpicture}[scale = 1.4]
    \node[wvert,label={right}:$\starw_{b_i}$] (w11) at (0,-2) {};
    \node[wvert,label={above}:$\starw_{b_{i+1}}$] (w22) at (-2,0) {};

    \coordinate[wvert] (w3) at (65:1.8){};
    \coordinate[wvert] (w4) at (25:1.8){};
    \node[] at (45:1.2) {$\ddots$};

    \node[bvert, label=below right:$\starb_{f_{i,i+1}}$] (b12) at (0,0) {};
    \draw[] (b12) edge node[sloped, above]{1} (w3)
                 edge node[sloped, below]{1} (w4);

    \path[font=\scriptsize]
      (w11) edge node [sloped, above] {$\lambda_i$} (b12)
      (w22) edge node [sloped, below] {$-(1+\lambda_{i+1})$} (b12);

    \node (fw) at (-1,-1) {$\starf_w$};
  \end{tikzpicture}
  \caption*{(B)}
  \label{fig:section_hexagon:B}
\end{subfigure}

\caption{Local configuration near the corner between $w$ and $f_{i,i+1}$ along with edge weights and vertex/face labels in $(\vrc,\relation)$ on $\pb$ (A) and in the corresponding VRC $(\starvrc,\starrelation)$ on $\stargraph$ (B).}
\label{fig:section_hexagon}
\end{figure}

The goal of this section is to show that the notion of a section ties together the projective and affine cluster structures introduced in Section~\ref{sec:clusterstructures}: the affine cluster structure of a TCD map is obtained as the projective cluster structure of its section.

\begin{theorem}\label{th:affprojcluster}
	Let $\tcd$ be a TCD and let $\tcdmap: \tcdm \rightarrow \CP^\dimension$ with $\dimension \geq 2$ be a 1-generic TCD map, and let $\hyperplane$ be a hyperplane that is generic with respect to $\tcdmap$. Then
    \[
    		\pro(\sigma_\hyperplane(\tcdmap)) = \aff_\hyperplane(\tcdmap).\qedhere
    \]
\end{theorem}
\begin{proof}
     Let $\projectivequiver$ and $\affinequiver$ denote the quivers of $\pro(\sigma_\hyperplane(\tcdmap))$ and $\aff_\hyperplane(\tcdmap)$ respectively. We first check that the quivers are the same.

    By Proposition~\ref{lem:star_graph_properties} (1), the white vertices of $\pb$ correspond to the faces of $\sigma(\pb)$ and boundary white vertices correspond to boundary faces, so both the mutable and the frozen vertices of the two quivers are the same. 
    
    We now show that they also have the same arrows. Since splits and contractions preserve the projective quiver, $\stargraph$ has the same projective quiver as $\sigma(\pb)$. Let $\starf_w$ denote the face of $\stargraph$ corresponding to a white vertex $w$ of $\pb$. There is an arrow in $\projectivequiver$ from $\starf_{w}$ to $\starf_{w'}$ if and only if there is an edge $\starb_f \starw_b$ in $\stargraph$ separating $\starf_{w}$ and $\starf_{w'}$ such that $\starf_{w}$ is on the right-hand side of the edge when it is oriented from $\starb_f$ to $\starw_b$, which is true if and only if $w,b,w'$ appear in counterclockwise cyclic order around $f$ in $\pbint$. There is a corresponding arrow in $\affinequiver$ from $w$ to $w'$ if $b$ is an internal black vertex of $\pbint$ by Step~(2) of Definition~\ref{def:affinecluster} and if $b$ is a boundary black vertex of $\pbint$ by Step~(3) of Definition~\ref{def:affinecluster}. When $f$ has degree $4$, cancellation of oriented $2$-cycles occurs in both quivers.

 In the special case of Step~(5) of Definition~\ref{def:star_graph}, 
 $\projectivequiver$ acquires a loop at the glued vertex which gets deleted in that step, while $\affinequiver$ has no arrow.

    Finally, we show that the mutable cluster variables in $\pro(\sigma_\hyperplane(\tcdmap))$ and $\aff_\hyperplane(\tcdmap)$ are identical. This can be verified locally at each internal white vertex $w$ of $\pb$.

    Let $m$ be the degree of $w$, and let $b_1,\dots,b_m$ be the black vertices adjacent to $w$ in counterclockwise cyclic order. For each $i$, let $w_i$ and $w_i'$ denote the two other white vertices adjacent to $b_i$. By applying a projective transformation, we may assume that $\hyperplane$ is the hyperplane at infinity. Let $(\vrc,\relation)$ be a VRC of $\tcdmap$ in affine gauge. Let $f_{i,i+1}$ be the face of $\pb$ incident to $w$ and located between $b_i$ and $b_{i+1}$.
    
    Since projective cluster variables are preserved under splits and contractions at degree-$2$ black vertices, we compute them in $\stargraph$. By applying gauge transformations, we may choose the edge weights of $(\vrc, \relation)$ as in Figure~\ref{fig:section_hexagon}~(A), corresponding to the relations
    \[
       \relation(b_i): \quad \vrc(w) + \lambda_i \vrc(w_i) - (1 + \lambda_i)\vrc(w'_i) =0.
    \]
    Note that all the following vectors are parallel to each other and each of them is a lift of $\sigma_\hyperplane(\tcdmap)(b_i)$:
    \[
        \vrc(w) - \vrc(w_i), \qquad \vrc(w) - \vrc(w_i'), \qquad \vrc(w_i)-\vrc(w_i').
    \]
    We choose a VRC lift $(\starvrc,\starrelation)$ on $\stargraph$ associated to the section $\sigma_{\hyperplane}(\tcdmap)$ by contraction of degree-$2$ white vertices as in Figure~\ref{fig:vrc_moves} as follows:
    \[
        \starvrc(\starw_{b_i}) = \vrc(w_i) - \vrc(w'_i), \qquad
        \starvrc(\starw_{b_{i+1}}) = \vrc(w_{i+1}) - \vrc(w'_{i+1}),
        \]
    while for $b \neq b_{i}, b_{i+1}$ on the boundary of $f_{i,i+1}$ with $w' \to b \to w''$ in counterclockwise order around $f_{i,i+1}$, 
    \[
        \starvrc(\starw_b) = \vrc(w'')-\vrc(w').
    \]
    The resulting relation $\starrelation(\starb_{f_{i,i+1}})$ in the section is obtained by first subtracting consecutive relations
    \[
         \relation(b_i) - \relation(b_{i+1}): \quad 
        \lambda_i(\vrc(w_i) - \vrc(w'_i))  - (1 + \lambda_{i+1})(\vrc(w_{i+1}) - \vrc(w'_{i+1}))+ (\vrc(w_{i+1}) - \vrc(w'_i))=0,
    \]
    and then rewriting 
    \[
        \vrc(w_{i+1}) - \vrc(w'_i) = \sum_{b \neq b_i,b_{i+1}} \starvrc(\starw_b),\qquad \text{$b \neq b_{i}, b_{i+1}$ on the boundary of $f_{i,i+1}$},
    \]
    to get 
    \[
        \starrelation(\starb_{f_{i,i+1}}): \quad 
        \lambda_i\starvrc(\starw_{b_i}) - (1 + \lambda_{i+1})\starvrc(\starw_{b_{i+1}}) + \sum_{b \neq b_i,b_{i+1}} \starvrc(\starw_b) =0.
    \]

    This gives the edge weights shown in Figure~\ref{fig:section_hexagon}~(B).
   
    The contribution of $f_{i,i+1}$ to the star-ratio $Y_w$ is \(-\frac{1+\lambda_i}{\lambda_{i+1}}\), and its contribution to the multi-ratio $X_{\starf_w}$ is \(-\frac{1+\lambda_{i+1}}{\lambda_i}\). Taking the product over all $i = 1, \dots, m$ and multiplying by $(-1)^{m+1}$, we get that $Y_w = X_{\starf_w}$.
\end{proof}

Since the section of a TCD map is again a TCD map, one may also consider iterated sections of a TCD map. Moreover, because the rank of $\sigma(\tcdmap)$ is at most the rank of $\tcdmap$ minus one, one can take at most $\maxrank(\tcdmap)$ iterated sections. A priori, this would suggest that up to $2\maxrank(\tcdmap)$ distinct cluster structures can be associated to a TCD map. However, Theorem~\ref{th:affprojcluster} implies that there are in fact at most $\maxrank(\tcdmap)+1$ distinct cluster structures.

\begin{remark}
    Note that, under our standing minimality assumption, we have defined sections only for TCD maps on minimal TCDs. Nevertheless, the same definition applies to non-minimal TCDs as well. However, in that case, some information may be lost when taking sections but Theorem~\ref{th:affprojcluster} continues to hold away from local violations of minimality.
\end{remark}

\section*{Acknowledgements}

We thank Pavlo Pylyavskyy, Wolfgang Schief, Boris Springborn, Ananth Sridhar and Yuri Suris for several fruitful discussions with some of us at various stages of the completion of this project. Part of this paper was written while TG was visiting the Indian Statistical Institute, Bangalore. NA was supported by the Deutsche Forschungsgemeinschaft (DFG) Collaborative Research Center TRR 109 “Discretization in Geometry and Dynamics”. NA and SR were partially supported by the Agence Nationale de la Recherche, Grant Number ANR-18-CE40-0033 (ANR DIMERS). SR was also partially supported by the CNRS grant Tremplin@INP, which funded a visit of NA to Paris-Saclay.

\bibliographystyle{alpha}
\bibliography{references}

\end{document}